\documentclass{amsart}

\usepackage{amssymb}

\newtheorem{theorem}{Theorem}[section]
\newtheorem{lemma}[theorem]{Lemma}

\newtheorem{claim}[theorem]{Claim}

\newtheorem{example}[theorem]{Example}

\numberwithin{equation}{section}

\newcommand{\re}{\text{Re\,}}
\newcommand{\im}{\text{Im\,}}
\newcommand{\Ci}{\text{Ci\,}}
\newcommand{\Si}{\text{Si\,}}
\newcommand{\sgn}{\text{sgn\,}}

\begin{document}
\setcounter{page}{1}
\title[Determinantal processes with number variance saturation]
{Determinantal processes with number variance saturation}
\author[K.~Johansson]{Kurt Johansson}

\address{
Department of Mathematics,
Royal Institute of Technology,
S-100 44 Stockholm, Sweden}

\email{kurtj@math.kth.se}

\begin{abstract}
Consider Dyson's Hermitian Brownian motion model
after a finite time $S$, where the process is started at $N$ equidistant
points on the real line. These $N$ points after time $S$ form a
determinantal process and has a limit as $N\to\infty$. This limting
determinantal proceess has the interesting feature that it shows 
number variance
saturation. The variance of the number of particles in an interval converges
to a limiting value as the length of the interval goes to infinity. 
Number variance saturation is also seen for example in the zeros of the
Riemann $\zeta$-function, \cite{Odl}, \cite{Be}. The process can also
be constructed using non-intersecting paths and we consider several variants
of this construction. One construction leads to a model which shows a
transition from a non-universal behaviour with number variance saturation
to a universal sine-kernel behaviour as we go up the line.

\end{abstract}

\maketitle

\section{Introduction}
The Bohigas-Gianonni-Schmidt conjecture, \cite{BGS}, says that the spectrum
$\{E_j\}_{j\ge 1}$, $E_j\to\infty$ as $j\to\infty$, of a quantum system whose 
classical dynamics is fully chaotic, has random matrix statistics in the large
energy limit, $E_j\to\infty$. For finite $E_j$ the specrum has non-universal 
features depending on the particular system, but as we go higher up in the 
spectrum the statistical properties become more and more like those from the
universal point processes obtained from random matrix theory. The zeros of
Riemann's $\zeta$-function, $\{E_j\}_{j\ge 1}$, $0<E_1<E_2<\dots$
(assuming the Riemann hypothesis) show a similar behaviour and has been 
popular as a model system in quantum chaos, \cite{BeKe}, since there are many
analogies. The number of zeros $\le E$, denoted by $\mathcal{N}(E)$, is 
approximately $\frac E{2\pi}\log\frac E{2\pi}-\frac E{2\pi}$. If we unfold 
the zeros by letting $x_j=\mathcal{N}(E_j)$, so that the mean spacing 
becomes 1, it is conjectured by Montgomery, \cite{Mo}, 
and tested numerically by
Odlyzko, \cite{Odl}, 
that the statistics of the $x_j$:s as $j\to\infty$, is like
the statistics of a determinantal point process with correlation functions
\begin{equation}\label{1.1}
\det\left(\frac{\sin\pi(x_i-x_j)}{\pi(x_i-x_j)}\right)_{i,j=1}^m,
\end{equation}
$m\ge 1$. Weak forms of this have been proved in \cite{He}, \cite{RuSa}.

This determinantal point process is obtained as a scaling limit of GUE or 
the classical compact groups, e.g. $U(n)$, and also as the universal scaling
limit of many other hermitian random matrix ensembles. If we count the number 
of particles (eigenvalues) in an interval of length $L$ in a process on 
$\mathbb{R}$ with correlation functions (\ref{1.1}) we get a random variable 
with variance, the {\it number variance}, $\sim\frac 1{\pi^2}\log L$ as
$L\to\infty$ (an exact formula for finite $L$ is given in (\ref{2.34'})
below). A feature in the quantum chaos model and for the Riemann zeros is 
{\it number variance saturation}, \cite{Be}, \cite{Se}. 
If we consider the Riemann zeros and intervals of 
length $L$ at height $E$, where $L\ll E$, and compute the variance by 
considering many disjoint intervals of length $L$, the dependence on $L$
is such that for small $L<\log\frac E{2\pi}$ it behaves like
$\frac 1{\pi^2}\log L$ but as $L$ grows it saturates, actually oscillates 
around an average value which is
approximately $\frac 1{\pi^2}\log(\log \frac E{2\pi})$,
see the work of Berry, \cite{Be}, for interesting precise predictions.
Hence the sine kernel 
determinantal point process is only a good model in a restricted range which 
becomes longer as we go up the line. The question that we address in this 
paper is whether it is possible to construct a determinantal process which
shows number variance saturation? Can we construct a determinantal process on
$[0,\infty)$ which shows a transition from a non-universal
regime to a universal regime described by (\ref{1.1})  as we go further 
and further away
from the origin? It is not possible to get number variance saturation 
with a translation invariant kernel, like in (\ref{1.1}), since the sine
kernel is the kernel with the slowest growth of the number variance among
all translation invariant kernels which define a determinantal point
process, \cite{So}. 

In this paper we will construct models having these properties by suitable 
scaling limits of determinantal process defined using 
non-intersecting Brownian motions. These models will not be translation 
invariant. We can restore translation invariance by averaging, but then we 
will no longer have a determinantal point process. One of the kernels 
obtained is given approximately by
\begin{equation}
\frac{\sin\pi(x-y)}{\pi(x-y)}+\frac{d\cos\pi(x+y)+(y-x)\sin\pi(y+x)}
{\pi(d^2+(y-x)^2)},
\end{equation}
where $d>0$ is a parameter, see (\ref{2.26}) below.
(The exact kernel has corrections of order $\exp(-Cd)$). This model will have
a number variance with saturation level $\sim\frac 1{\pi^2}\log(2\pi d)$.

There are other connections between $L$-functions and random matrix theory.
Katz and Sarnak, \cite{KaSa}, study low-lying zeros of families of L-functions
and connect their statistical behaviour with that obtained for the 
eigenvalues of random matrices from the compact classical groups with
respect to Haar measure. This leads to a classification of the L-functions
into different symmetry classes. Three different laws for the distribution of
the lowest zero are obtained. We will see below that these three laws can
also be obtained from the non-intersecting Brownian motions by choosing
different boundary conditions, see also \cite{Gra}. 
Another recent development is the study
of characteristic polynomials of matrices from the classical groups
which have been used to model L-functions, and led to interesting conjectures 
for their moments, see \cite{KeSn}. See also \cite{CoDi} for a discussion
of linear statistics of zeros.

In this paper we will have nothing to say about quantum chaos or L-funtions. 
The above discussion only serves as a background and a motivation for 
dicussing the probabilistic models that we will introduce. 
For a discussion of bounded variance in another context see \cite{AGL}.

\section{Models and results}
\subsection{The model}

A point process on $B\subseteq \mathbb{R}$ with correlation functions 
$\rho_n(x_1,\dots, x_n)$, $n\ge 1$, \cite{So}, has determinantal
correlation functions if there is a function $K:B\times B\to\mathbb{R}$,
the {\it correlation kernel}, such that
\begin{equation}\label{2.0}
\rho_n(x_1,\dots, x_n)=\det(K(x_i,x_j))_{i,j=1}^n,
\end{equation}
$n\ge 1$. The interpretation of $\rho_n$ is that $\rho_n(x_1,\dots,x_n)
dx_1\dots dx_n$ is the probability of finding particles in infinitesimal
intervals $dx_1,\dots, dx_n$ around $x_1,\dots x_n$. In particular
$\rho_1(x)=K(x,x)$ is the local density at $x$.

Below we will construct kernels $K$ by taking appropriate limits of other 
kernels and it is natural to ask if there is a determinantal point process 
whose correlation kernel is $K$. This can be answered using the following
theorem.

\begin{theorem}\label{thm2.0}
Let $\rho_{k,N}(x_1,\dots,x_k)$, $k\ge 1$, be the correlation functions of a
determinantal point process on $I\subseteq \mathbb{R}$ with continuous
correlation kernel $K_N(x,y)$, $N\ge 1$. Assume that $K_N(x,y)\to K(x,y)$
uniformly on compact subsets of $I^2$. Then there is a point process on 
$I$ with correlation functions
\begin{equation}\label{2.0'}
\rho_k(x_1,\dots, x_k)=\det(K(x_i.x_j))_{i,j=1}^k,
\end{equation}
$k\ge 1$.
\end{theorem}

The theorem will be proved at the end of sect. 3.

Let $\phi_i(t)$, $\psi_i(t)$, $i\ge 1$ be functions in $L^2(X,\mu)$. Then,
\cite{Bo}, \cite{TW}, 
\begin{equation}\label{2.1}
u_N(x)d^N\mu(x)=\frac 1Z\det(\phi_i(x_j))_{i,j=1}^N
\det(\psi_i(x_j))_{i,j=1}^Nd^N\mu(x)
\end{equation}
defines a measure on $X^N$ with determinantal correlation functions. 
If we have $u_N(x)\ge 0$ for all $x\in X^N$ and
\begin{equation}\label{2.2}
Z=\int_{X^N}\det(\phi_i(x_j))_{i,j=1}^N
\det(\psi_i(x_j))_{i,j=1}^Nd^N\mu(x)>0,
\end{equation}
we get a probability measure. We can 
think of a symmetric probability 
measure on $X^N$ as a point process on $X$ with exactly $N$ particles.
In this paper we will have $X=\mathbb{R}$ or $[0,\infty)$ and $\mu$ will be
Lebesgue measure. The correlation kernel is given by
\begin{equation}\label{2.3}
K_N(x,y)=\sum_{i,j=1}^N \psi_i(x)(A^{-1})_{ij}\phi_j(y),
\end{equation}
where $A=(\int_{X}\phi_i(x)\psi_j(x)d\mu (x))_{i,j=1}^N$. Note that
$\det A=Z>0$.

One natural way to obtain probability measures of the form (\ref{2.1}) is 
from non-intersecting paths using the Karlin-McGregor theorem, \cite{KaMc}. 
Consider
$N$ one-dimensional Brownian motions started at $y_1<\dots<y_N$ at time 0
and conditioned to stop at $z_1<\dots<z_N$ at time $S+T$ and not to 
intersect (coincide) in the whole time interval $[0,S+T]$. The induced
measure on the positions $x_1,\dots, x_N$ at time $S$ is then
\begin{equation}\label{2.4}
p_{N,S,T}(x)=\frac 1Z\det(p_S(y_i,x_j))_{i,j=1}^N\det(p_T(x_i,z_j))_{i,j=1}^N
\end{equation}
where 
\begin{equation}
p_t(x,y)=\frac 1{\sqrt{2\pi t}}e^{-(x-y)^2/2t}
\notag
\end{equation}
is the transition kernel for one-dimensional Brownian motion. From the results
discussed above it follows that (\ref{2.4}) defines a point process on
$\mathbb{R}$ with determinantal correlation functions. The correlation
kernel is given by
\begin{equation}\label{2.5}
K_{N,S,T}(u,v)=\sum_{k,j=1}^Np_S(y_k,u)(A^{-1})_{kj}p_T(v,z_j),
\end{equation}
where
\begin{equation}\label{2.6}
A=\det(p_{S+T}(y_i,z_j))_{i,j=1}^N.
\end{equation}
In the limit $T\to\infty$ this model converges to Dyson's Brownian motion 
model, \cite{Dy}, with $\beta=2$, and was considered in \cite{Jo}, see also
\cite{Guhr}.

We can also consider the same type of measure but on $[0,\infty)$ and with 
appropriate boundary conditions at the origin. We will consider reflecting
or absorbing boundary conditions, where we have the transition kernels,
\cite{Brei},
\begin{equation}\label{2.7}
p_t^{\text{re}}(x,y)=p_t(x,y)+p_t(x,-y)
\end{equation}
and
\begin{equation}\label{2.8}
p_t^{\text{ab}}(x,y)=p_t(x,y)-p_t(x,-y)
\end{equation}
respectively. We simply replace $p_S$, $p_T$ in (\ref{2.4}) with 
$p_S^{\text{re}}$, $p_T^{\text{re}}$ or 
$p_S^{\text{ab}}$, $p_T^{\text{ab}}$. In these cases we have initial points
$0<y_1<\dots<y_N$ and final points $0<z_1<\dots<z_N$. We will be interested
in these models as $N\to\infty$ with fixed $S,T$ or with $S$ fixed and
$T\to\infty$.

\subsection{Correlation kernels}

We want to obtain useful expressions for the correlation kernel (\ref{2.5})
with equidistant final positions, compare \cite{Jo}, proposition 2.3.

\begin{theorem}\label{thm2.1}
Let $\Gamma_L:\mathbb{R}\ni t\to L+it$, $L\in\mathbb{R}$, and let $\gamma$ be
a simple closed curve that surrounds $y_1<\dots <y_{2n-1}$; $L$ is so large 
that $\Gamma_L$ and $\gamma$ do not intersect. Set $z_j=a(j-n)$ for some
$a>0$, $1\le j\le 2n-1$. Consider the model (\ref{2.4}) with $(y_i)_{i=0}
^{2n}$ as initial conditions, and $(z_i)_{i=0}^{2n}$ as final points, and
with no boundary. Then
\begin{align}\label{2.15}
&K_{2n+1,S,T}(u,v)=\frac{ae^{-(u^2+v^2)/2T}}{(2\pi i)^2S(S+T)}\int_{\Gamma_L}
dw\int_{\gamma}dz \frac 1{e^{a(w-z)/(T+S)}-1}
\\
&e^{-\frac{T}{2S(T+S)}[(w-\frac{T+S}Tv)^2+(z-\frac{T+S}Tu)^2]
+\frac {an}{T+S}(z-w)} \prod_{j=0}^{2n-1}\frac{e^{aw/(T+S)}-e^{ay_j/(T+S)}}
{e^{az/(T+S)}-e^{ay_j/(T+S)}}.
\notag\end{align}
\end{theorem}

The theorem will be proved in sect. 3.

The $T\to\infty$ limit of this formula appears in \cite{Jo}, compare theorem
\ref{thm2.2} below. There are analogues of the formula (\ref{2.15}) for the
absorbing and reflecting cases. We have not been able to write down a
useful formula for the case of general final positions. The expression in
(\ref{2.15}) 
is more useful computationally than (\ref{2.3}) but still
rather complicated. We will obtain simpler formulas in
certain special cases. First we will give a double contour integral formula
for the case $T=\infty$ in the absorbing and reflecting cases. We will also
consider the $N\to\infty$ formula in the absorbing case. Then we will
specialize to the case when the initial points 
are also equidistant and $T=\infty$ or
$T=S$. In these last two cases we can obtain very nice formulas that are
not in terms of contour integrals. The next theorem gives the analogue of
proposition 2.3 in \cite{Jo} in the absorbing and reflecting cases.

\begin{theorem}\label{thm2.2}
Let $\Gamma_L$ be as in theorem \ref{thm2.1} and assume that $\gamma$ 
surrounds $y_1,\dots,y_N$ and does not intersect $\Gamma_L$. Set $z_j=j-1$
and assume $0<y_1<\dots<y_N$. Then, uniformly for $(u,v)$ in a compact set in
$[0,\infty)^2$,
\begin{align}\label{2.16}
&\lim_{T\to\infty}K_{N,S,T}^{\text{ab}}(u,v)=K_{N,S}^{\text{ab}}(u,v)\\
&\doteq \frac 1{(2\pi i)^2S}\int_{\Gamma_L}dw\int_{\gamma} dz
e^{(w-v)^2/2S}(e^{-(z-u)^2/2S}-e^{-(z+u)^2/2S})\frac {2w}{w^2-z^2}
\prod_{j=1}^N \frac{w^2-y_j^2}{z^2-y_j^2}
\notag\end{align}
and
\begin{align}\label{2.17}
&\lim_{T\to\infty}K_{N,S,T}^{\text{re}}(u,v)=K_{N,S}^{\text{re}}(u,v)\\
&\doteq \frac 1{(2\pi i)^2S}\int_{\Gamma_L}dw\int_{\gamma} dz
e^{(w-v)^2/2S}(e^{-(z-u)^2/2S}+e^{-(z+u)^2/2S})\frac {2z}{w^2-z^2}
\prod_{j=1}^N \frac{w^2-y_j^2}{z^2-y_j^2}.
\notag\end{align}
\end{theorem}

The theorem will be proved in sect. 3. 

We will also write down a contour integral formula for the $N\to\infty$
limit of $K_{N,S}^{\text{ab}}(u,v)$ under an assumption on the $y_j$:s.
(We could write a similar formula in the reflecting case, but we will omit 
it.)

\begin{theorem}\label{thm2.3}
Let $\Gamma_L$ be as above, $L$ arbitrary, 
and $\gamma_M$ the two lines $\mathbb{R}\ni t\to \mp t\pm iM$ with $M>0$.
Let $0<y_1<y_2<\dots$ and assume that $\sum_{j=1}^\infty 1/y_j^2<\infty$.
Define
\begin{equation}\label{2.18}
F(z)=\prod_{j=1}^\infty \left(1-\frac{z^2}{y_j^2}\right),
\end{equation}
which converges uniformly on all compact subsets of $\mathbb{C}$. Set
\begin{equation}\label{2.19}
K^\ast_{S,1}(u,v)=\frac 1{(2\pi i)^2S}\int_{\Gamma_L} dw\int_{\gamma_M} dz
e^{(w-v)^2/2S-(z-u)^2/2S}\frac 1{z-w}\frac{wF(w)}{zF(z)},
\end{equation}
\begin{equation}\label{2.20}
K^\ast_{S,2}(u,v)=\frac 1{2\pi i}\int_{L-Mi}^{L+Mi}e^{(w-v)^2/2S-
(w-u)^2/2S}dw,
\end{equation}
and
$K_S^\ast(u,v)=K^\ast_{S,1}(u,v)+K^\ast_{S,2}(u,v)$. Then uniformly on
compact subsets of $[0,\infty)^2$,
\begin{equation}\label{2.21}
\lim_{N\to\infty} K_{N,S}^{\text{ab}}(u,v)=K_{S}^{\text{ab}}(u,v)\doteq
K_S^\ast(u,v)-K_S^\ast(-u,v)
\end{equation}
\end{theorem}

The theorem will be proved in section 3. 

We now come to the case when the initial points $y_j$ are equidistant. The 
next theorem is what makes it possible to compute the number variance in this 
case. To compute the number variance using the double contour integrals
seems difficult.

\begin{theorem}\label{thm2.4}
Let $y_j=\Delta+a(j-n)$, $1\le j\le 2n-1$, $0\le \Delta<a$, $a>0$. Set
\begin{equation}\label{2.22}
d=\frac {2\pi S}{a^2}.
\end{equation}
Then, uniformly for $(u,v)$ in a compact subset of $\mathbb{R}^2$,
\begin{equation}\label{2.22'}
\lim_{N\to\infty} K_{N,S}(u,v)=K_S(u-\Delta,v-\Delta),
\end{equation}
where $K_S(u,v)=a^{-1}L_S(a^{-1}u,a^{-1}v)$, and
\begin{equation}\label{2.23}
L_S(x,y)=\frac 1{\pi}\re \sum_{n\in\mathbb{Z}} e^{-\pi dn(n-1)}
\frac {e^{\pi i (y+(2n-1)x)}}{nd+i(y-x)}.
\end{equation}
Furthermore, if $y_j=aj$, $j\ge 1$, then uniformly on compact subsets
of $[0,\infty)^2$,
\begin{equation}\label{2.24}
\lim_{N\to\infty} K^{\text{ab}}_{N,S}(u,v)=K^{\text{ab}}_S(u,v)\doteq
K_S(u,v)-K_S(-u,v),
\end{equation}
and
\begin{equation}\label{2.25}
\lim_{N\to\infty} K^{\text{re}}_{N,S}(u,v)=K^{\text{re}}_S(u,v)\doteq
K_S(u,v)+K_S(-u,v).
\end{equation}
\end{theorem}

The theorem will be proved in sect. 3.

The leading contribution to (\ref{2.23}) comes from the terms $n=0$ and $n=1$.
The other terms are exponentially small in $d$. The leading part is
\begin{equation}\label{2.26}
L_{S,\text{appr}}(x,y)=\frac {\sin \pi(x-y)}{\pi(x-y)}
+\frac {d\cos \pi(x+y)+(y-x)\sin \pi(y+x)}{\pi(d^2+(y-x)^2)},
\end{equation}
so we have the ordinary sine kernel plus a non-translation invariant term.
In particular for the density we have
\begin{equation}\label{2.26'}
K_{S,\text{appr}}(u,u)=\frac 1a \left(1+\frac{\cos(2\pi u/a)}{\pi d}\right),
\end{equation}
so we have an oscillating density reflecting the initial configuration.

It is also possible to give a more explicit formula in the case $T=S$.
We have the following theorem.
\begin{theorem}\label{thm2.5}
Let $\theta_i(x;\omega)$, $i=1,2,3,4$, be the Jacobi theta functions, see
(\ref{3.35}), 
let $y_j=\Delta+a(j-n)$, $1\le j\le 2n-1$, $0\le\Delta<a$, $a>0$,
and let $d$ be as in (\ref{2.22}). Then, uniformly on compact subsets
of $\mathbb{R}$,
\begin{equation}\label{2.27}
\lim_{N\to\infty} K_{N,S,S}(u,v)=K_{S,S}(u-\Delta,v-\Delta),
\end{equation}
where $K_{S,S}(u,v)=a^{-1}L_{S,S}(a^{-1}u,a^{-1}v)$ and 
\begin{align}\label{2.28}
L_{S,S}(u,v)=&\frac 1{\theta_2(0;id)\sqrt{\theta_3(0;id)\theta_4(0;id)}}
\\
&\times\left(\frac{\theta_3(u+v;2id)\theta_1(u-v;2id)}
{\sinh(\pi (u-v)/2d)}+
\frac{\theta_2(u+v;2id)\theta_4(u-v;2id)}{\cosh(\pi (u-v)/2d)}\right).
\notag\end{align}
\end{theorem}

If we neglect contributions which are exponentially small in $d$ the
leading part of (\ref{2.28}) is
\begin{equation}\label{2.29}
L_{S,S,\text{appr}}(u,v)=\frac{\sin\pi(u-v)}{2d\sinh(\pi (u-v)/2d)}
+\frac{\cos\pi(u+v)}{2d\cosh(\pi (u-v)/2d)}.
\end{equation}
Note that as $d\to\infty$ both $L_S$ and $L_{S,S}$ converge to the sine 
kernel. We see from (\ref{2.26}) and (\ref{2.29}) that $L_{S,S}$ decays much 
faster than $L_S$ at long distances.

We can also consider an averaged model by averaging over $\Delta$ in theorems
\ref{thm2.4} and \ref{thm2.5}. The averaged model has correlation functions
\begin{equation}\label{2.30}
\frac 1a \int_0^a \det(K(x_i-\Delta,x_j-\Delta))_{i,j=1}^m d\Delta,
\end{equation}
where $K$ is the appropriate kernel $K_S$ or $K_{S,S}$. This averaging will 
restore translation invariance. The density in the averaged process will be 
a constant equal to $1/a$. 
In particular we have
\begin{align}\label{2.31}
&\frac 1a\int_0^a K_S(x-\Delta,y-\Delta)K_S(y-\Delta,x-\Delta)d\Delta
\\&=
\frac{\sin^2\frac{\pi(x-y)}a}{\pi^2(x-y)^2}+
\frac{d^2-\left(\frac{x-y}a\right)^2}
{2\pi^2a^2(d^2+\left(\frac{x-y}a\right)^2)^2}.
\notag
\end{align}
plus terms exponentially small in $d$. 
If we instead consider $K_{S,S}$ we get
\begin{equation}\label{2.32}
\frac{\sin^2\frac{\pi (x-y)}a}{4a^2d^2\sinh^2\frac{\pi(x-y)}{2ad}}+
\frac 1{8a^2d^2\cosh^2\frac{\pi(x-y)}{2ad}}.
\end{equation}

\subsection{The number variance}

Let $I\subseteq\mathbb{R}$ be an interval and denote by $\#I$ the 
number of particles contained in $I$.
We are interested 
in the variance, $Var_K(\# I)$, of this random variable
in the determinantal point process with kernel $K$. The kernel
$K=K_{N,S,T}$ that we have considered above is a reproducing kernel, i.e.
\begin{equation}\label{2.33}
\int_{\mathbb{R}} K(x,y)K(y,z)dy=K(x,z).
\end{equation}
This is immeidiately clear from (\ref{2.5}) and (\ref{2.6}) and the same also
holds in the absorbing and reflecting cases with integration over $[0,\infty)$
instead, and the reproducing property is inherited by the limiting kernels 
obtained above. Using (\ref{2.33}) and the determinantal form of the
correlation functions it follows that
\begin{equation}\label{2.34}
Var_K(\# I)=\int_{I}dx\int_{I^c}dy K(x,y)K(y,x)
\end{equation}
The sine kernel $\frac{\sin\pi(x-y)/a}{\pi (x-y)}$ with density $1/a$ has
the number variance, $I=[0,L]$,
\begin{equation}\label{2.34'}
V_{\text{sinekernel}}(L)=\frac 1{\pi^2}\left[\log\frac{2\pi L}a +\gamma+1
+\frac{2\pi L}{a}(\frac {\pi}2 -\text{Si\,}(\frac{2\pi L}a))-
\cos\frac{2\pi L}a-\text{Ci\,}(\frac{2\pi L}a)\right]
\end{equation}
In the averaged models we get (we denote the averaging over $\Delta$
by $<\,>$),
\begin{equation}\label{2.35}
<\text{Var}_K(\# I)>=\int_Idx\int_{I^c} dy\frac 1a \int_0^a K(x-\Delta,
y-\Delta)K(y-\Delta,x-\Delta) d\Delta,
\end{equation}
where $K$ is $K_S$ or $K_{S,S}$. The formulas above for the correlation
functions and the formulas for the number variance are used to prove the
next theorems. We will only consider the contributions from the leading
parts of the kernels, (\ref{2.26}) and (\ref{2.29}). Also, we will not use 
the reflecting and absorbing kernels. If the intervals are high up, 
$I=[R,R+L]$ with $R$ large, then the contribution from $K_S(-u,v)$ in 
(\ref{2.24}) and (\ref{2.25}) will be small (like $O(1/R)$). 

\begin{theorem}\label{thm2.6}
Consider the kernel $K_S$ of theorem \ref{thm2.4}. Define $\theta$ and $\phi$
by $R/a-[R/a]=\theta/\pi$ and $L/a-[L/a]=\phi/\pi$. The contribution 
to $\text{Var}_{K_S}(\#[R,R+L])$ coming from the leading part of the
kernel (\ref{2.26}) is, $A=\pi L/a$,
\begin{align}\label{2.36}
&\frac 1{\pi^2}\left(1+\frac{\cos 2(\theta+\phi)+\cos 2\theta}{\pi d}\right)
\left(\log\frac {2\pi Ad}{\sqrt{A^2+\pi^2d^2}}+\gamma-\text{Ci\,}(2A)\right)
\\&+
\frac 1{\pi^2}\left(1+2A\left(\frac {\pi}2-\text{Si\,}(2A)\right)
-\cos 2A\right)
\notag\\
&+\frac 1{2\pi^3d}\left\{(\cos 2(\theta+\phi)+\cos 2\theta)(h_3(2i\pi d)
-h_3(2A+2i\pi d))\right.
\notag\\&\left.+
(\sin 2(\theta+\phi)-\sin 2\theta)(h_1(2A+2i\pi d)+\pi
-2\text{Si\,}(2A))\right\}
\notag\\
&+\frac{\sin 2\phi}{8\pi^3 d}\left\{h_2(2A+2i\pi d)
(\sin 4(\theta+\phi)-\sin 4\theta)-h_4(2A+2i\pi d)
(\cos 4(\theta+\phi)+\cos 4\theta)\right\},
\notag\end{align}
where $h_1(z)=2\re f(z)$, $h_2(z)=2\im f(z)$, $h_3(z)=2\re g(z)$,
$h_4(z)=2\im g(z)$ and
\begin{equation}\label{2.37}
f(z)=\int_0^\infty \frac{\sin t}{t+z} dt\quad ,\quad
g(z)=\int_0^\infty \frac{\cos t}{t+z} dt,
\end{equation}
for $\re z>0$ or $\re z=0$, $\im z\neq 0$. These functions have the 
asymptotics $f(z)=1/z-2/z^3+O(1/z^4)$, $g(z)=1/z^2+ O(1/z^4)$ as
$z\to\infty$.
\end{theorem}

The theorem will be proved in sect. 4.

If we average the expression for the variance over $\theta$ or equivalently
consider the averaged model, we get that the contribution to the number
variance from the leading part is
\begin{equation}\label{2.38}
\frac 1{\pi^2}\left[\log\frac{2\pi Ld/a}{\sqrt{L^2/a^2+d^2}}+\gamma+1
-\text{Ci\,}(\frac{2\pi L}a)+\frac{2\pi L}a (\frac {\pi}2 -
\text{Si\,}(\frac{2\pi L}a)-\cos\frac{2\pi L}a)\right].
\end{equation}
Apart from the logarithmic term we have exactly the same formula as for
the sine kernel, (\ref{2.34'}). It is not difficult to obtain (\ref{2.38})
using the averaged correlation functions, (\ref{2.31}). The proof is then
analogous to that in the sine kernel case, see sect. 4. 
The proof of (\ref{2.35})
is a lengthy but rather straightforward computation. 

Note that when $L$ is small compared to $d$ the leading term in (\ref{2.38})
is $\frac 1{\pi^2}\log(2\pi L/a)$, which is what we have for the sine kernel.
When $L\to\infty$ the expression (\ref{2.38}) converges to 
\begin{equation}\label{2.39}
\frac 1{\pi^2}(\log(2\pi d)+\gamma+1)
\end{equation}
so the number variance saturates. Note that the saturation level does
not depend directly on the mean spacing. If we rescale the model, see below, 
we get the same saturation level.

We can also compute the number variance for the (S,S)-model where we have
the kernel $K_{S,S}$. In this case we will only consider the averaged model.

\begin{theorem}\label{thm2.7}
Denote by $V_d(L)$ the leading part of the 
averaged number variance of an interval of length
$L$ for the model with kernel $K_{S,S}$, i.e. (\ref{2.29}), where $d$ 
is given by (\ref{2.22}). Then
\begin{align}\label{2.40}
V_d(L)&=\frac 2{\pi^2}\int_{0}^{\pi L/a}\left(\frac{u/2\pi d}{\sinh(u/2\pi d)}
\right)^2\frac{\sin^2 u}{u} du+\frac {2L}{\pi a}\int_{\pi L/a}^\infty
\left(\frac{u/2\pi d}{\sinh(u/2\pi d)}
\right)^2\frac{\sin^2 u}{u^2} du
\\
&+\frac 1{\pi^2}\int_0^\infty\frac{\min(u,L/2ad)}{\cosh^2 u} du.
\notag\end{align}
Also,
\begin{equation}\label{2.41}
\lim_{d\to\infty}\frac 1{\frac 1{\pi^2}\log d} \lim_{L\to\infty} V_d(L) =1.
\end{equation}
\end{theorem}
If we compare (\ref{2.40}) with the integrals which lead to (\ref{2.34'}),
we see that there is a truncation effect which depends on $d$ and which 
is responsible for the saturation. The limit (\ref{2.41}) shows that the
saturation level is similar to (\ref{2.39}) for large $d$.

As mentioned in the introduction the unitary group $U(n)$ has been used in 
\cite{KeSn} to model the $\zeta$-function at height $T$, where 
$n=\log(T/2\pi)$. This $n$ is obtained by equating the mean spacing in $U(n)$,
$2\pi /n$, with the mean spacing of the zeros at height $T$, which is 
$2\pi(\log(T/2\pi))^{-1}$. Note that the eigenvalues of a random matrix
from $U(n)$ with respect to the Haar measure also show a kind of number 
variance saturation. The variance for the number of eigenvalues in an interval
on $\mathbb{T}$ of length $a$, $0\le a\le 2\pi$ is given by
\begin{equation}\label{2.41'}
\text{Var\,}_{U(n)}(a)=\frac {na}{2\pi}-\frac{na^2}{4\pi^2}-
\frac 2{\pi^2}\sum_{k=1}^{n-1}\frac{n-k}{k^2}\sin^2\frac {ka}2.
\end{equation}
This increases as a function of $a$ for $0\le a\le \pi$ and then decreases
symmetrically. we have a maximum variance when $a=\pi$, i.e when we
have a half circle. This maximum variance is
\begin{equation}\label{2.42}
\frac 1{\pi^2}(\log (2n)+\gamma+1)+O(\frac 1n),
\end{equation}
which is analogous to (\ref{2.39}) if we set $\pi d=n$. Note that the number
variance for the Riemann zeros at height $T$ saturates at the mean level
$\frac 1{\pi^2}\log(\log(T/2\pi))+\text{const}$, \cite{Be}, so equating
the saturation levels (disregarding constant terms) leads to 
$n=\log(T/2\pi)$ again. This may be a more natural argument in a sense
since it is not changed under rescaling.

\subsection{Approximation}

Consider the situation in theorem \ref{thm2.3}, 
the absorbing $(S,\infty)$-model,
where the $y_j$'s are not equally spaced but are given by $y_j=F^{-1}(j)$ 
for some increasing function $F$. If $F$ is nice and do not vary too quickly,
the $y_j$'s will be almost equally spaced for long stretches of $j$, and 
hence we expect that the kernel (\ref{2.21}) should be well approximated by
(\ref{2.24}) in a region where the average spacing is $a$. We will not 
attempt to make this clear in the greatest possible generality. Our goal is
an approximation theorem valid for a certain class of functions $F$. 
Denote the correlation kernel $K_S^{\text{ab}}(u,v)$ with initial points
$\underline{y}=(y_j)_{j=1}^\infty$ by $K_S^{\text{ab}}(\underline{y},u,v)$
to indicate the dependence on $\underline{y}$. We will prove the following 
approximation theorem:

\begin{theorem}\label{thm2.8}
Assume that $F:[0,\infty)\to[0,\infty)$ is a $C^2$-function that satisfies
\newline
(i) $F(x)\le Cx^{1+\delta}$, for all $x\ge 0$, for some $\delta\in (0,1)$
and some constant $C$,
\newline
(ii) $F'(0)=0$ and $F'(x)>0$ for $x>0$,
\newline
(iii) $F''$ is decreasing.
\newline
Fix $\alpha$ (large) and $\epsilon >0$ so that $1+\delta+\epsilon<2$.
Also, fix $S\in[F'(F(\alpha), 1]$. Define $\underline{y}=(y_j)_{j=1}^\infty$
by $y_j=F^{-1}(j)$, $j\ge 1$ and $y_{-j}=-y_j$ for $j\ge 0$. For
$m\ge 1$ set
\begin{equation}\label{2.43}
\xi_m=\sum_{j=1}^\infty\frac{2y_m-y_{m-j}-y_{m+j}}{(y_m-y_{m-j})
(y_{m+j}-y_m)},
\end{equation}
$\lambda_m^{-1}=F'(y_m)$ and $\eta_m=F''(y_m)$. Also let $\zeta_m=
y_m-\xi_mS$. There is an $m=m(\alpha)\ge 1$ such that
$|\zeta_m-\alpha|\le \lambda_m$. We have roughly $m(\alpha)\sim F(\alpha)$.
Set $\lambda(\alpha)=\lambda_{m(\alpha)}\sim (F'(F^{-1}(\alpha))^{-1}$,
$\eta(\alpha)=\eta_{m(\alpha)}\sim F''(F^{-1}(\alpha))$ and
\begin{equation}\label{2.44}
T_0(\alpha)=\min(\frac 1{4\sqrt{\eta(\alpha)}},m(\alpha)^{(1+\delta)^2
/2(1-\delta)}).
\end{equation}
Define $\tilde{\underline{y}}=\{\tilde{y}_j\}_{j=1}^\infty$, by 
$\tilde{y}_j=\lambda(\alpha)j$, $j\ge 1$. Assume that $0<T\le T_0(\alpha)$.
There are constants $c_0$, $C$, which depend on $F$ and $\epsilon$, but not
$\alpha$, such that
\begin{align}\label{2.45}
&\left|e^{-\xi(\alpha)(u-v)}K_S^{\text{ab}}(\underline{y};u,v)-
K_S^{\text{ab}}(\tilde{\underline{y}};u-\zeta(\alpha),v-\zeta(\alpha))\right|
\\
&\le C[\lambda(\alpha)S^{-3/2}e^{-R^2/8S}+(T^2+R^2)m(\alpha)^{-(1-\delta)^2/
(1+\delta)}S^{-1}],
\notag\end{align}
for all $u,v\in[\alpha-T,\alpha+T]$ provided $R^2$ lies in the interval
\begin{equation}\label{2.46}
\left[c_0^2\max\left(\left(\frac S{\lambda(\alpha)}\right)^{2/(1-\epsilon)},
T^{1+\delta+\epsilon} S,
\frac{T^{1+\epsilon}S}{\lambda(\alpha)}\right)\,,\,
m(\alpha)^{(1-\delta)^2/(1+\delta)}\right].
\end{equation}
\end{theorem}

Note that $K_S^{\text{ab}}(\tilde{\underline{y}};x,y)$ is given by 
(\ref{2.24}) with $a=\lambda(\alpha)$, so provided the right hand side of
(\ref{2.45}) is small we have an approximation with a kernel having
equally spaced initial points. The factor $\exp(-\xi(\alpha)(u-v))$
in front of $K_S^{\text{ab}}(\underline{y};u,v)$ in (\ref{2.45}) does not
affect the correlation functions corresponding to this kernel since it
cancels in the determinant. Let us consider two examples of theorem 2.8.

\begin{example}\label{ex2.9}\rm
Let $F(x)=x^{1+\delta}$ with $0<\delta<1$ and fix $\epsilon >0$ small.
Then $\lambda(\alpha)^{-1}\sim (1+\delta)\alpha^{\delta/(1+\delta)}$,
$\eta(\alpha)\sim\delta(1+\delta)\alpha^{-(1-\delta)/(1+\delta)}$
and $m(\alpha)\sim\alpha^{1+\delta}$ as $\alpha\to\infty$. Choose $S=1$. If
$0<T\le \alpha^{(1-\delta)^2/2}$, then
\begin{equation}\label{2.47}
\left|e^{-\xi(\alpha)(u-v)}K_1^{\text{ab}}(\underline{y};u,v)-
K_1^{\text{ab}}(\tilde{\underline{y}};u-\zeta(\alpha),v-\zeta(\alpha))\right|
\le C\frac{T^2}{\alpha^{(1-\delta)^2-2\delta/[(1-\epsilon)(1+\delta)]}} ,
\end{equation}
for all $u,v\in [\alpha-T,\alpha+T]$. We see that this is only interesting if 
$\delta<2-\sqrt{3}$, otherwise the right hand side of (\ref{2.47}) does not
go to zero as $\alpha\to\infty$ unless we also let $T\to 0$. It also
follows that 
\begin{equation}\label{2.48}
\lim_{\alpha\to\infty}
e^{-\xi(\alpha)(x-y)/\lambda(\alpha)} K_1(\underline{y};\alpha+\frac 
x{\lambda(\alpha)},\alpha+\frac y{\lambda(\alpha)})=\frac {\sin\pi (x-y)}
{\pi(x-y)}
\end{equation}
uniformly for $x,y$ in a compact set, so as we go up
the line we see the sine kernel process. Note that when $x$ and $y$
belong to a compact set we can take $T=C\lambda(\alpha)$ and the right
hand side of (\ref{2.47}) goes to zero as $\alpha\to\infty$ for $\delta<1$
provided we choose $\epsilon$ small enough. Thus we can extend (\ref{2.48})
to all $0<\delta<1$
We see a transition from a 
non-universal regime for small $\alpha$ to the universal sine kernel 
regime for large
$\alpha$.\it
\end{example}

\begin{example}\label{ex2.10} \rm
Let $F(x)=\frac x{2\pi}\log\frac x{2\pi}-
\frac x{2\pi}$. This $F$ does not satisfy all the conditions in theorem 
\ref{thm2.8} but we can modify it for small $x$ so that it does without
changing the $y_j$:s. Then $\lambda(\alpha)\frac 1{2\pi}\log
\frac {\alpha}{2\pi}\to 1$, $2\pi\alpha\eta(\alpha)\to 1$ and $m(\alpha)/
F(\alpha)\to 1$ as $\alpha\to\infty$. Fix $\epsilon_0>0$ small and let
$\delta=\epsilon$ where $(1-\epsilon)^2/(1+\epsilon)-2\epsilon=1-\epsilon_0$.
Pick $R=c_0T^{1/2+\epsilon}\alpha^\epsilon$. The theorem then shows that
if $1\le T\le \alpha^{1/2-\epsilon}$, then ($S=1$)
\begin{equation}\label{2.49}
\left|e^{-\xi(u-v)}K_1(\underline{y};u,v)-K_1(\tilde{\underline{y}};u,v)
\right|\le CT^2\alpha^{-1+\epsilon}
\end{equation}
for $u,v\in[\alpha-T,\alpha+T]$. We see that the parameter $d$ for the
aproximate kernel is
\begin{equation}
d(\alpha)=\frac{2\pi}{\lambda(\alpha)^2}\sim\frac {8\pi^3}{(\log(\alpha/2\pi)
)^2}.
\notag
\end{equation}
The saturation level is the (disregarding constants)
\begin{equation}\label{2.50}
\sim\frac 1{\pi^2}\log d(\alpha)\sim \frac 2{\pi^2}\log(\log\frac
{\alpha}{2\pi}).
\end{equation}
This differs by a factor 2 from what we would like to have if we want to 
model the $\zeta$-function. Note that the mean spacing at height $\alpha$ 
is the same as for the $\zeta$-function. However the saturation level does 
not change in the equidistant case if we rescale the point process linearly
(the distance $a$ between the $y_j$:s is not changed, it is the constructed
point process that is rescaled). Hence the mean spacing and the saturation
level are independent. We can make a better model of the $\zeta$-function
by picking a suitable $F$, see below, and then make a non-linear 
rescaling. Below we will construct a model for the unfolded zeros $s_j=
\mathcal{N}(E_j)$.

Consider a determinantal point process on $[0,\infty)$ with correlation kernel
$K$ and let $G(x)$ be a strictly increasing $C^1$-function. Then
\begin{equation}\label{2.51}
\tilde{K}(x,y)=K(G(x),G(y))\sqrt{G'(x)G'(y)}
\end{equation}
defines a new, rescaled point process. The density for the new process is
$\rho(G(x))G'(x)$ if $\rho(x)$ is the density for the original process.

Let $y_j=F^{-1}(j)$, $j\ge 1$, where $F$ satisfies the conditions of theorem
\ref{thm2.8}. The density at $x$ is then $F'(x)$ and if we rescale with $G$
we get $F'(G(x))G'(x)=\frac d{dx}F(G(x))$. Hence to get a constant density
we should pick $G(x)=F^{-1}(x)$. Choose
\begin{equation}\label{2.52}
F(x)=1+\int_{2\pi e}^x(\log\frac t{2\pi})^{1/2}dt
\end{equation}
for $x\ge 2\pi e$ and define it for $0\le x\le 2\pi e$ so that the conditions
in theorem \ref{thm2.8} are satisfied. Consider the model corresponding to
this $F$ and rescale it using $G=F^{-1}$ as in (\ref{2.51}). We will say 
``approximately'' below without beeing too precise. The estimates involved
can be made precise with a little effort. By theorem \ref{thm2.8} at
height $T$, the kernel will be approximately ($S=1$)
\begin{equation}\label{2.53}
K_1(\tilde{\underline{y}};F^{-1}(x)-F^{-1}(T), F^{-1}(y)-F^{-1}(T))
\sqrt{(F^{-1})'(x)(F^{-1})'(y)}
\end{equation}
for $x,y$ in a neighbourhood of $T$, where $\tilde{y}_j=\lambda(F^{-1}(T))j$,
$j\ge 1$, $\lambda(\alpha)^{-1}\sim F'(\alpha)$. (We have approximated
$\zeta(\alpha)$ with $\alpha$.) Now, for $x$ close to $T$ we have
$F^{-1}(x)-F^{-1}(T)\approx (x-T)(F^{-1})'(T)$. Set $b=(F^{-1})'(T)\approx
\lambda(T)$. By (\ref{2.53}) and (\ref{2.26}) at height $T$ (large) the kernel
will approximately equal
\begin{equation}\label{2.54}
L_{S,\text{appr}}(x-T,y-T)=\frac {\sin\pi(x-y)}{\pi(x-y)}+
\frac{d\cos\pi(x+y-2T)+(y-x)\sin\pi(x+y-2T)}{\pi(d^2+(y-x)^2)},
\end{equation}
where $d=2\pi/\lambda(T)^2$. We have $\lambda(T)\sim(\log\frac{T}{2\pi})^
{-1/2}$. As $T\to\infty$ , $d\to\infty$ and we see that the kernel in
(\ref{2.54}) converges to the sine kernel. The point process we have 
constructed thus has the correct universal asymptotics as we go up
the line, and it is non-universal for small $T$, since it depends on the
particular $F$ we have chosen. The saturation level for the kernel in
(\ref{2.54}) is
\begin{equation}\label{2.55}
\frac 1{\pi^2}\log (2\pi d)\sim\frac 1{\pi^2}\log (\log\frac T{2\pi}),
\end{equation}
when $T$ is large, which is exactly what we would like to have at height $T$.
We have thus constructed a determinantal point process in $[0,\infty)$
which has, in some aspects, similar behaviour to the unfolded zeros of the
$\zeta$-function.\it
\end{example}

\subsection{Correlation kernels close to the origin}
Consider the kernels $K_S$ on $\mathbb{R}$, and $K_S^{\text{ab}}$ and
$K_S^{\text{re}}$ on $[0,\infty)$ in theorem \ref{thm2.4}. If we are 
interested in say the distribution of the first particle to the right of
the origin we have to compute the probability of having no particle in
$[0,\xi]$. Let $X$ be the position of the first particle to the right
of the origin. Then, if the correlation kernel is $K$,
\begin{align}\label{2.56}
\mathbb{P}[X\le\xi]&=1-\mathbb{P}[\text{no particle in $[0,\xi]$}]
\\
&=1-\det(I-K)_{L^2([0,\xi])},
\notag
\end{align}
where the second equality is a standard result for determinantal point
processes. If we
vary $a$ and $S$ in such a way that $S/a^2\to\infty$, it follows from 
theorem \ref{thm2.4} that
\begin{equation}\label{2.57}
\lim_{S/a^2\to\infty} aK_S(au,av)=\frac{\sin\pi(u-v)}{\pi(u-v)},
\end{equation}
\begin{equation}\label{2.58}
\lim_{S/a^2\to\infty} aK_S^{\text{ab}}(au,av)=\frac{\sin\pi(u-v)}{\pi(u-v)}
-\frac{\sin\pi(u+v)}{\pi(u+v)}
\end{equation}
and
\begin{equation}\label{2.59}
\lim_{S/a^2\to\infty} aK_S^{\text{re}}(au,av)=\frac{\sin\pi(u-v)}{\pi(u-v)}
+\frac{\sin\pi(u+v)}{\pi(u+v)}
\end{equation}
These kernels can also be obtained from the classical compact groups and have
been used by Katz and Sarnak to model the lowest zeros in families of
L-functions, see \cite{KaSa}. The above results show that they can
also be otained in a natural way from non-intersecting paths with appropriate
boundary conditions.

The kernels in the right hand side of (\ref{2.58}) and (\ref{2.59}) are 
directly related to special instances of the Bessel kernel,
\begin{equation}\label{2.60}
B_\nu(x,y)=\frac{x^{1/2}J_{\nu+1}(x^{1/2})J_{\nu}(y^{1/2})
-J_{\nu}(x^{1/2})y^{1/2}J_{\nu+1}(y^{1/2})}{2(x-y)},
\end{equation}
where $J_\nu(x)$ is the ordinary Bessel function. In fact a simple computation
shows that if we define the rescaled Bessel kernel $\tilde{B}_\nu$ by
\begin{equation}\label{2.61}
\tilde{B}_\nu(x,y)=\sqrt{2\pi^2 x2\pi^2y}B_\nu(\pi^2x^2,\pi^2y^2),
\end{equation}
then
\begin{equation}\label{2.62}
\tilde{B}_{\pm 1/2}(x,y)=
\frac{\sin\pi(u-v)}{\pi(u-v)}\mp \frac{\sin\pi(u+v)}{\pi(u+v)}.
\end{equation}
When $\nu$ is an integer the kernel $B_\nu$ appears in the scaling limit
around the smallest eigenvalue in LUE, the Laguerre Unitary Ensemble. 
If $x_1,\dots,x_N$ are the eigenvalues
of $M^\ast M$, where $M$ is a $(\nu+N)\times N$ complex matrix with 
independent standard complex Gaussian elements, $N(0,1/2)+iN(0,1/2)$, then
$x_1,\dots,x_N$ is a finite determinantal point process with correlation
kernel $K_N^\nu$ and
\begin{equation}\label{2.63}
\lim_{N\to\infty}\frac 1{4N}K^\nu_N(\frac x{4N},\frac y{4N})=B_\nu(x,y).
\end{equation}
This interpretation does not work for $\nu=\pm 1/2$.

\section{Computation of the correlation functions}
In this section we will use the formula (\ref{2.5}) to compute the 
correlation functions. If $A$ is a matrix and $b$ a column vector we
will denote by $(A|b)_k$ the matrix where column $k$ in $A$ is replaced
by $b$. Let
\begin{equation}
p=(p_T(v,z_0)\dots p_T(v,z_N))^T,
\notag
\end{equation}
where $z_j=a(j-n)$, $0\le j\le N=2n$. By (\ref{2.5}) and Kramers rule we have
\begin{equation}\label{3.1}
K_{2n,S,T}(\underline{y};u,v)=\sum_{k=0}^{2n}p_S(y_k,u)\frac{\det(A|p)_k}
{\det A},
\end{equation}
where $A=(p_{S+T}(y_i,z_k))_{j,k=0}^{2n}$. If $C_M$ is the contour $t\to
t+iM$, $M\in\mathbb{R}$, $t\in\mathbb{R}$, we have
\begin{equation}\label{3.2}
p_T(v,z_i)=\sqrt{\frac{T+S}{2\pi T}}e^{-v^2/2T}\int_{C_M}e^{-\tau^2/2}
p_{T+S}(\tilde{y}_k(\tau),z_i)e^{\tilde{y}_k(\tau)^2/2(T+S)}d\tau,
\end{equation}
where $\tilde{y}_k=\tilde{y}_k(\tau)=v(T+S)/T+i\tau\sqrt{S(T+S)/T}$;
set also $\tilde{y}_i=y_i$ if $i\neq k$. Then,
\begin{equation}\label{3.3}
\frac{\det(A|p)_k}{\det A}=\sqrt{\frac{T+S}{2\pi T}}
e^{-v^2/2T}\int_{C_M}e^{-\tau^2/2+\tilde{y}_k(\tau)^2/2(T+S)}
\frac{\det(p_{T+S}(\tilde{y}_j,z_i))_{i,j=0}^{2n}}
{\det(p_{T+S}(y_j,z_i))_{i,j=0}^{2n}} d\tau.
\end{equation}
Since $z_j=a(j-n)$ the determinants in the quotient in the right hand side of
(\ref{3.3}) can be computed using Vandermonde's determinant and we find
\begin{equation}
e^{(y_k^2-\tilde{y}_k^2)/2(T+S)+an(y_k-\tilde{y}_k)/2(T+S)}
\prod_{j\neq k}\frac{e^{a\tilde{y}_k/(T+S)}-e^{ay_j/(T+S)}}
{e^{ay_k/(T+S)}-e^{ay_j/(T+S)}}.
\notag
\end{equation}
Inserting this into (\ref{3.2}) and making the change of variables 
$w=\tilde{y}_k(\tau)$ we obtain
\begin{align}\label{3.4}
&\frac{\det(A|p)_k}{\det A}=\frac 1{i\sqrt{2\pi S}}e^{-v^2/2T+y_k^2/2(T+S)+
ny_ka/(T+S)}
\\
&\times \int_{\Gamma_L}dw e^{T(w-(T+S)v/T)^2/2S(T+S)-nwa/(T+S)}
\prod_{j\neq k}\frac{e^{aw/(T+S)}-e^{ay_j/(T+S)}}
{e^{ay_k/(T+S)}-e^{ay_j/(T+S)}}.
\notag\end{align}
We can now use the expression in (\ref{3.4}) and insert it into (3.1) to
get
\begin{align}\label{3.5}
&K_{2n,S,T}(\underline{y};u,v)=\frac 1{2\pi S}e^{-v^2/2T+u^2/2T}
\sum_{k=0}^{2n}e^{-T(y_k-(T+S)u/T)^2/2S(T+S)}
\\
&\times \int_{\Gamma_L}dw e^{T(w-(T+S)v/T)^2/2S(T+S)}
e^{na(y_k-w)/(T+S)}
\prod_{j\neq k}\frac{e^{aw/(T+S)}-e^{ay_j/(T+S)}}
{e^{ay_k/(T+S)}-e^{ay_j/(T+S)}}.
\notag\end{align}
This is the basic formula from which the others will be derived. It is now
straightforward to prove theorem \ref{thm2.1}.

\begin{proof}
{\it of theorem \ref{thm2.1}}. If we apply the residue theorem in (\ref{2.15})
we get the expression (\ref{3.5}). Since $\Gamma_L$ and $\gamma$ do not
intersect the $w=z$ singularity does not contribute.
\end{proof}

We turn next to theorem \ref{thm2.2}.
\begin{proof}
{\it of theorem \ref{thm2.2}}. Consider first the absorbing case. We use the 
formulas (\ref{3.1}) and (\ref{3.3}) but with $p_t^{\text{ab}}$ instead
of $p_t$. The evaluation of the two determinants can now be done using the
following Vandermonde type identity
\begin{equation}\label{3.6}
\det(e^{(2i-1)x_j}-e^{-(2i-1)x_j})_{i,j=1}^N=
\prod_{i=1}^N(e^{x_i}-e^{-x_i})\prod_{1\le i<j\le N}((e^{x_j}-e^{-x_j})^2-
(e^{x_i}-e^{-x_i})^2).
\end{equation}
Using this identity we find
\begin{align}\label{3.7}
\frac{\det(A|p^{\text{ab}})_k}{\det A}&=\sqrt{\frac{T+S}{2\pi T}}
e^{\frac{y_k^2}{2(T+S)}-\frac{v^2}{2T}}\int_{C_M} e^{-\frac{\tau^2}2}
\frac{e^{\frac vT+i\sqrt{\frac{S}{T(T+S)}}}-
e^{-\frac vT-i\sqrt{\frac{S}{T(T+S)}}}}{e^{\frac{y_k}{T+S}}-
e^{-\frac{y_k}{T+S}}}
\\
&\times\prod_{j\neq k}\frac{\left(e^{\frac vT+i\sqrt{\frac{S}{T(T+S)}}}-
e^{-\frac vT-i\sqrt{\frac{S}{T(T+S)}}}\right)^2-
\left(e^{\frac{y_j}{T+S}}-
e^{-\frac{y_j}{T+S}}\right)^2}
{\left(e^{\frac{y_k}{T+S}}-e^{-\frac{y_k}{T+S}}\right)^2
-\left(e^{\frac{y_j}{T+S}}-e^{-\frac{y_j}{T+S}}\right)^2}.
\notag
\end{align}
In this expression we can expand the exponentials and take the $T\to\infty$
limit. We find that
\begin{equation}\label{3.8}
\lim_{T\to\infty}\frac{\det(A|p^{\text{ab}})_k}{\det A}=
\frac 1{i\sqrt{2\pi S}}\int_{\Gamma_L} e^{(w-v)^2/2S}\frac w{y_k}
\prod_{j\neq k}\frac{w^2-y_j^2}{y_k^2-y_j^2} dw
\end{equation}
uniformly for $(u,v)$ in a compact set. Here we have changed integration 
variable by putting $w=v+i\tau\sqrt{S}$. Hence
\begin{align}\label{3.9}
&\lim_{T\to\infty}K^{\text{ab}}_{N,S,T}(u,v)
\\
&=\frac 1{2\pi iS}\sum_{k=1}^N\left(e^{-(y_k-u)^2/2S}-e^{-(y_k+u)^2/2S}
\right)\int_{\Gamma_L} e^{(w-v)^2/2S}\frac w{y_k}
\prod_{j\neq k}\frac{w^2-y_j^2}{y_k^2-y_j^2} dw.
\notag\end{align}
That this expression equals the expression in (\ref{2.16}) follows from
the residue theorem.

The proof of (\ref{2.17}) is completely analogous. Instead of (\ref{3.6}) 
we use the identity
\begin{equation}\label{3.10}
\det \left(e^{(j-1)x_i}-e^{-(j-1)x_i}\right)_{i,j=1}^N
=\prod_{1\le i<j\le N}((e^{x_j}+e^{-x_j})-(e^{x_i}-e^{-x_i}).
\end{equation}
\end{proof}

\begin{proof} ({\it of theorem \ref{2.3}}). Set
\begin{equation}\label{3.11}
F_N(z)=\prod_{j=1}^N\left(1-\frac{z^2}{y_j^2}\right),
\end{equation}
so that, by (\ref{2.16}),
\begin{equation}\label{3.12}
K^{\text{ab}}_{N,S}(u,v)=\frac 1{(2\pi i)^2S}\int_{\Gamma_0} dw
\int_{\gamma^+_{M,r}} dze^{\frac{(w-v)^2}{2S}}\left(e^{\frac{(z-u)^2}{2S}}
-e^{\frac{(z+u)^2}{2S}}\right)\frac {2w}{w^2-z^2}\frac{F_N(w)}{F_N(z)}.
\end{equation}
Here $\gamma^+_{M,r}$ is the curve given by
$t\to-t\pm im$, $t\ge r$ and $s\to r\mp is$, $-M\le s\le M$, where
$0<r<y_1$. The fact that $F_N(z)\to F(z)$ uniformly on compact subsets
of $\mathbb{C}$, together with estimates like (\ref{5.9}) and (\ref{5.10})
below, which can be used to restrict the $z-$ and $w-$ integrations, shows
that
\begin{align}
\lim_{N\to\infty}K^{\text{ab}}_{N,S}(u,v)&=\frac 1{(2\pi i)^2S}
\int_{\Gamma_0} dw
\int_{\gamma^+_{M,r}} dze^{\frac{(w-v)^2}{2S}}
\\
&\times\left(e^{\frac{(z-u)^2}{2S}}
-e^{\frac{(z+u)^2}{2S}}\right)\frac {2w}{w^2-z^2}\frac{F(w)}{F(z)}
\doteq K_S^{\text{ab}}(u,v).
\notag
\end{align}
uniformly for $(u,v)$ in a compact set. Let $\gamma_{M,r}^-$ be the contour
which is the image of $\gamma^+_{M,r}$ under $z\to -z$. Using
\begin{equation}
\frac {2w}{w^2-z^2}=\frac{w}{z}\left(\frac 1{w-z}-\frac 1{w+z}\right)
\notag
\end{equation}
we see that $K_S^{\text{ab}}(u,v)=K_S^\ast(u,v)-K_S^\ast(-u,v)$, where
\begin{equation}
K_S^\ast(u,v)=\frac 1{(2\pi i)^2S}\int_{\Gamma_0} dw\int_{\gamma_{M,r}}dz
e^{\frac{(w-v)^2}{2S}-\frac{(z-u)^2}{2S}}\frac 1{z-w}\frac{wF(w)}{zF(z)},
\notag
\end{equation}
$\gamma_{M,r}=\gamma^+_{M,r}+\gamma^-_{M,r}$. If we let $r\to 0+$ we pick
up a contribution from the pole $z=w$. This leads to $K^\ast_{S}=
K^\ast_{S,1}+K^\ast_{S,2}$ as stated in the theorem.
\end{proof}

Next we consider theorem \ref{thm2.4}.

\begin{proof} {\it of theorem \ref{thm2.4}}. 
It suffices to consider the case $\Delta=0$, otherwise we replace $(u,v)$ 
with  $(u-\Delta,v-\Delta)$. From the proof of proposition 2.3 in \cite{Jo},
the $T\to\infty$ limit of (\ref{3.5}) is
\begin{align}\label{3.13}
K_{N,S}(u,v)&=\lim_{T\to\infty} K_{N,S,T}(u,v)
\\
&=\frac 1{2\pi i S}\sum_{k=-n}^n e^{-\frac {(ak-u)^2}{2S}}
\int_{\Gamma_L} dw e^{\frac {(w-v)^2}{2S}}\prod_{j=-n, j\neq k}^n
\frac{w-aj}{ak-aj}.
\notag\end{align}
Now,
\begin{equation}\label{3.14}
\prod_{j=-n, j\neq k}^n\frac{w-aj}{ak-aj}=\frac{(-1)^kw}{w-ak}
\prod_{j=1}^n\left(1-\frac{w^2}{(aj)^2}\right)\prod_{j=0}^{k-1}
\frac{n-j}{n+1+j}.
\end{equation}
It follows from (\ref{3.13}) and (\ref{3.14}) that
\begin{align}\label{3.15}
&\lim_{N\to\infty}K_{N,S}(u,v)=K_S(u,v)
\\
&=\frac 1{2\pi i S}\sum_{k\in\mathbb{Z}} e^{-\frac {(ak-u)^2}{2S}}
\int_{\Gamma_L} dw e^{\frac {(w-v)^2}{2S}}
\prod_{j=1}^\infty\left(1-\frac{w^2}{(aj)^2}\right)
\notag\end{align}
uniformly for $(u,v)$ in a compact set. To prove the convergence we need
some estimates so that we can cut off the $k$-summation and the 
$w$-integration. We omit the details. Now,
\begin{equation}\label{3.16}
w\prod_{j=1}^\infty\left(1-\frac{w^2}{(aj)^2}\right)=
\frac a{\pi}\sin\frac{\pi}a w.
\end{equation}
Hence,
\begin{equation}\label{3.17}
K_S(u,v)=\frac{a}{2\pi^2iS}\sum_{k\in\mathbb{Z}}\int_{\Gamma_L}
e^{\frac{(w-v)^2-(ak-u)^2}{2S}}\frac{(-1)^k}{w-ak}\sin\frac{\pi w}a dw,
\end{equation}
where $L$ is arbitrary. Replace $w$ by $w+v-u+ak$ in the integral in
(\ref{3.17}) and use Cauchy's theorem to get
\begin{align}\label{3.18}
K_S(u,v)&=\frac{a}{2\pi^2iS}\sum_{k=0}^\infty 
\int_{\Gamma_{-1}}e^{\frac{w^2-2uw}{2S}}
e^{\frac {kaw}S}\frac{\sin\frac {\pi}{a}(w+v-u)}{w+v-u} dw
\\
&+\frac{a}{2\pi^2iS}\sum_{k=-\infty}^{-1} 
\int_{\Gamma_{1}}e^{\frac{w^2-2uw}{2S}}
e^{\frac {kaw}S}\frac{\sin\frac {\pi}{a}(w+v-u)}{w+v-u} dw
\notag\\
&=\frac{a}{2\pi^2iS}\int_{\Gamma}\frac{e^{(w^2-2uwD)/2S}}{e^{aw/S}-1}
\frac{\sin\frac {\pi}{a}(w+v-u)}{w+v-u} dw,
\notag\end{align}
where $\Gamma=\Gamma_{-1}+\Gamma_1$. The function $e^{aw/S}-1$ has simple
zeros at $w=2\pi Sni/a$, $n\in\mathbb{Z}$ and hence by the residue theorem
applied to the last integral in (\ref{3.18}),
\begin{equation}\label{3.19}
K_S(u,v)=\frac 1{\pi}\sum_{n\in\mathbb{Z}} e^{-2\pi^2n^2S/a^2-2\pi uni/a}
\frac{\sin\frac {\pi}a\left(\frac{2\pi nS}a i+v-u\right)}
{\frac{2\pi nS}a i+v-u}.
\end{equation}
If we set $L(x,y)=aK(ax,ay)$ and define $d$ as in (\ref{2.22}) we find
\begin{align}
L(x,y)&=\frac 1{\pi}\sum_{n\in\mathbb{Z}} e^{-\pi dn^2-2\pi nix}
\frac{\sin\pi(ndi+y-x)}{ndi+y-x}
\\
&=\frac 1{2\pi}\sum_{n\in\mathbb{Z}}\frac{e^{-\pi dn^2-\pi nd}}
{-nd+i(y-x)} e^{\pi i(y-(2n+1)x)}
\notag\\
&-\frac 1{2\pi}\sum_{n\in\mathbb{Z}}\frac{e^{-\pi dn^2+\pi nd}}
{-nd+i(y-x)} e^{\pi i(y+(2n-1)x)}.
\notag
\end{align}
If we change $n$ to $-n$ in the first sum we get (\ref{2.23}).

In the absorbing case we have by (\ref{3.9}), (\ref{3.11}) and the identity
\begin{equation}
\frac w{y_k}\prod_{j\neq k}\frac{w^2-y_j^2}{y_k^2-y_j^2}=
\left(\frac 1{w-y_k}-\frac 1{w+y_k}\right)\frac{F_N(w)}{F_N'(y_k)},
\notag
\end{equation}
that
\begin{align}\label{3.20}
\lim_{N\to\infty} K_{N,S}^{\text{ab}}(u,v)&=\frac a{2\pi^2iS}
\sum_{k=1}^\infty\int_{\Gamma_L}\left(e^{-\frac{(ak-u)^2}{2S}}-
e^{-\frac{(ak+u)^2}{2S}}\right)
\\
&\times 
\left(\frac 1{w-y_k}-\frac 1{w+y_k}\right)(-1)^k\sin\frac {\pi}a wdw.
\notag
\end{align}
If we now note that
\begin{align}
&\sum_{k=1}^\infty (-1)^k\left(e^{-\frac{(ak-u)^2}{2S}}-
e^{-\frac{(ak+u)^2}{2S}}\right)\left(\frac 1{w-y_k}-\frac 1{w+y_k}\right)
\\
&=\sum_{k\in\mathbb{Z}}e^{-\frac{(ak-u)^2}{2S}}\frac{(-1)^k}{w-ak}
-\sum_{k\in\mathbb{Z}}e^{-\frac{(ak+u)^2}{2S}}\frac{(-1)^k}{w-ak}
\notag\end{align}
then (\ref{3.17}) and (\ref{3.20}) give (\ref{2.24}). The proof of 
(\ref{2.25}) is completely analogous.
\end{proof}

The proof of theorem \ref{thm2.5} is a similar but somewhat more complicated 
computation where we have to use $\theta$-function identities.

\begin{proof} {\it of theorem \ref{thm2.5}}. Our starting point is the 
formula (\ref{3.5}) with $T=S$. We can assume that $\Delta=0$ so that 
$y_j=a(j-n)$, $0\le j\le 2n$. Set
\begin{equation}\label{3.21}
F_n(z;q)=\prod_{j=0}^n(1-q^j z)\prod_{j=1}^n(1-q^j/z),
\end{equation}
where
\begin{equation}\label{3.22}
q=e^{-a^2/2S}=e^{-\pi/d}.
\end{equation}
Then (\ref{3.5}) and a short computation gives
\begin{align}\label{3.23}
K_{2n+1,S,S}(u,v)&=\frac{ae^{-(v^2-u^2)/2S}}{2(2\pi i)^2 S^2}
\int_{\Gamma_L} dw\int_\gamma dz\frac 1{e^{a(w-z)/2S}-1}
\\
&\times e^{[(w-2v)^2-(z-2u)^2]/2S}\frac{F_n(e^{aw/2S};q)}{F_n(e^{az/2S};q)},
\notag\end{align}
where $\gamma$ surrounds $y_j=a(j-n)$, $0\le j\le 2n$ and does not
interesect $\Gamma_L$. Let $\gamma_M$ be $t\to\mp t\pm iM$. Replacing 
$\gamma$ by $\gamma_M$ in (\ref{3.23}) we pick up a contribution from the 
pole $z=w$, when $|\im w|\le M$. Hence
\begin{align}\label{3.24}
K_{2n+1,S,S}(u,v)&=\frac{ae^{-(v^2-u^2)/2S}}{2(2\pi i)^2 S^2}
\int_{\Gamma_L} dw\int_{\gamma_M} dz\frac 1{e^{a(w-z)/2S}-1}
\\
&\times e^{[(w-2v)^2-(z-2u)^2]/4S}\frac{F_n(e^{aw/2S};q)}{F_n(e^{az/2S};q)}
\notag\\
&-\frac{e^{(v^2-u^2)/2S}}{2\pi i}\int_{L-iM}^{L+iM} 
e^{[(w-2v)^2-(z-2u)^2]/4S}dw.
\notag
\end{align}
In this expression we can control the $n\to\infty$ limit (using some estimates
of $F_n$, compare (\ref{5.9}), (\ref{5.10})). Set
\begin{equation}\label{3.25}
F(z;q)=\prod_{j=0}^\infty (1-q^j z)\prod_{j=1}^\infty (1-q^j/z).
\end{equation}
It follows from (\ref{3.24}) that
\begin{align}\label{3.26}
K_{S,S}(u,v)=\lim_{N\to\infty}K_{2n+1,S,S}(u,v)
&=\frac{ae^{-(v^2-u^2)/2S}}{2(2\pi i)^2 S^2}
\int_{\Gamma_L} dw\int_{\gamma_M} dz\frac 1{e^{a(w-z)/2S}-1}
\\
&\times e^{[(w-2v)^2-(z-2u)^2]/4S}\frac{F(e^{aw/2S};q)}{F(e^{az/2S};q)}
\notag\\
&-\frac{e^{(v^2-u^2)/2S}}{2\pi i}\int_{L-iM}^{L+iM} 
e^{[(w-2v)^2-(z-2u)^2]/4S}dw.
\notag
\end{align}
We now compute the $z$-integral in (\ref{3.26}) using the residue theorem.
Apart from the pole $z=w$ if $|\im w|\le M$ we have simple poles at $z=ak$, 
$k\in\mathbb{Z}$. We obtain
\begin{align}\label{3.27}
K_{S,S}(u,v)&=
\frac{e^{-(v^2-u^2)/2S}}{2\pi iS}\sum_{k\in\mathbb{Z}}
\int_{\Gamma_L} dw e^{[(w-2v)^2-(z-2u)^2]/4S}
\\
&\times 
\frac 1{q^{-w/a+k}}\frac{F(q^{-w/a})}{q^{-k}F'(q^{-k})}.
\notag\end{align}
In the integral in (\ref{3.27}) we make the change of variables $w\to a(w+k)$
and use Cauchy's theorem. A computation shows that
\begin{equation}\label{3.28}
F(q^{-k}q^{-w})=\frac{(-1)^k}{q^{k(k+1)/2}}q^{-kw}F(q^{-w})
\end{equation}
\begin{equation}\label{3.29}
q^{-k}F'(q^{-k})=\frac{(-1)^{k-1}}{q^{k(k+1)/2}}
\prod_{n=1}^\infty (1-q^n)^2,
\end{equation}
which gives
\begin{equation}
\frac 1{q^{-w}-1}\frac{F(q^{-k}q^{-w})}{q^{-k}F'(q^{-k})}=
\prod_{n=1}^\infty (1-q^n)^{-2}\frac{q^{-kw}}{1-q^{-w}} F(q^{-w}).
\notag
\end{equation}
If we set 
\begin{equation}\label{3.30}
G(z;q)=\prod_{j=1}^\infty (1-q^j z)(1-q^j/z),
\end{equation}
we obtain
\begin{align}
K_{S,S}(u,v)&=\frac{ae^{-(v^2-u^2)/2S}}{2\pi iS}\sum_{k\in\mathbb{Z}}
\int_{\Gamma_L}dw e^{\frac 1{4S}[(aw+ak-2v)^2-(ak-2u)^2]}
\\
&\times e^{ka^2w/2S}\frac{G(q^{-w};q)}{G(1;q)}.
\notag\end{align}
Make the change of variables $w\to w+(v-u)/a$ and perform the $k$-summation
to get
\begin{equation}\label{3.31}
K_{S,S}(u,v)=-\frac{ae^{-(v^2+u^2)/2S}}{2\pi iS}\int_\Gamma
\frac{e^{(aw-u-v)^2/S}}{e^{-a^2w/S}-1}\frac{G(q^{-w+(u-v)/a};q)}
{G(1;q)} dw,
\end{equation}
where $\Gamma=\Gamma_1-\Gamma_{-1}$ as before. The integrand in (\ref{3.31})
has simple poles at $w=dni$, $n\in\mathbb{Z}$, and the residue theorem gives
\begin{equation}
K_{S,S}(u,v)=-\frac{e^{-(v^2+u^2)/2S}}{a}\sum_{n\in\mathbb{Z}}
e^{(adni-u-v)^2/S}\frac{G(q^{-dni+(u-v)/a};q)}{G(1;q)}.
\end{equation}
A computation leads to
\begin{align}\label{3.32}
L_{S,S}(u,v)&=aK_{S,S}(au,av)
\\
&=e^{-\pi(u-v)^2/2d}\sum_{n\in\mathbb{Z}} e^{-\pi dn^2/2-\pi ni(u+v)}
\frac{G((-1)^n e^{-\pi (u-v)/d};q)}{G(1;q)}.
\notag\end{align}
If we write $q=e^{\pi i\omega}$, $\omega=i/d$, $p=e^{\pi i x}$,
$C_0=\prod_{n=1}^\infty (1-q^{2n})$, it follows from the product 
representations of the $\theta$-functions that
\begin{align}\label{3.32'}
G(p^2;q)&=\frac 1{2C_0^2q^{1/4}}\frac{\theta_1(x;\omega)\theta_4(x;\omega)}
{\sin\pi x}
\\
G(-p^2;q)&=\frac 1{2C_0^2q^{1/4}}\frac{\theta_2(x;\omega)\theta_3(x;\omega)}
{\cos\pi x}
\notag\end{align}
Write $\omega'=-2/\omega=2id$.
If we insert the formulas (\ref{3.32'}) into (\ref{3.32}) we obtain
\begin{align}\label{3.33}
L_S(u,v)&=\frac{e^{-\frac{\pi (u-v)^2}{2d}}}{2C_0^2q^{1/4}G(1;q)}
\sum_{n\in \mathbb{Z}}\left\{e^{i\pi\omega' n^2-2\pi ni(u+v)}
\frac{\theta_1(\frac{\omega}2 (u-v);\omega)
\theta_4(\frac{\omega}2 (u-v);\omega)}{\sin\frac{\pi\omega}2 (u-v)}\right.
\\
&+\left.e^{i\pi\omega' (n-\frac 12)^2-\pi (2n-1)i(u+v)}
\frac{\theta_2(\frac{\omega}2 (u-v);\omega)
\theta_3(\frac{\omega}2 (u-v);\omega)}{\cos\frac{\pi\omega}2 (u-v)}\right\}
\notag\\
&=\frac{e^{-\frac{\pi (u-v)^2}{2d}}}{2C_0^2q^{1/4}G(1;q)}
\left\{\frac{\theta_3(u+v;\omega')\theta_1(\frac{\omega}2 (u-v);\omega)
\theta_4(\frac{\omega}2 (u-v);\omega)}{\sin\frac{\pi\omega}2 (u-v)}\right.
\notag\\
&\left.+\frac{\theta_2(u+v;\omega')\theta_2(\frac{\omega}2 (u-v);\omega)
\theta_3(\frac{\omega}2 (u-v);\omega)}{\cos\frac{\pi\omega}2 (u-v)}\right\},
\notag\end{align}
where we have used the series expansions of the $\theta$-functions, see 
(\ref{3.35}) below. 

We can simplify (\ref{3.33}) somewhat by using some $\theta$-function
identities. The first is Jacobi's transformations:

\begin{equation}
\theta_k(\omega x;\omega)=\frac{\theta_k(x;-1/\omega)}{\alpha_k
(-i\omega)^{1/2}e^{\pi i\omega x^2}},
\notag
\end{equation}
$k=1,2,3,4$, where $\alpha_1=-i$ and $\alpha_2=\alpha_3=\alpha_4=1$, 
and the formulas \cite{Law}, p. 17,
\begin{align}
&\theta_1(x;\omega)\theta_2(x;\omega)=\theta_1(2x;2\omega)\theta_4(0;2\omega),
\\
&\theta_3(x;\omega)\theta_4(x;\omega)=\theta_4(2x;2\omega)\theta_4(0;2\omega).
\notag
\end{align}
This leads to
\begin{equation}\label{3.34}
L_S(u,v)=\frac{i\theta_4(0;-2/\omega)}{2\omega C_0^2q^{1/4}G(1;q)}
\left\{\frac{\theta_3(u+v;\omega')\theta_1(u-v;\omega')}{-i\sin
\frac{\pi\omega}2 (u-v)}+
\frac{\theta_2(u+v;\omega')\theta_4(u-v;\omega')}{-i\cos
\frac{\pi\omega}2 (u-v)}\right\}.
\end{equation}
Now, by the product formulas for the values of the $\theta$-functions at
the origin we have
\begin{equation}
G(1;q)=\prod_{n=1}^\infty (1-q^n)^2=\frac{\theta_4(0;\omega)
\theta_1'(0,\omega)}{2\pi C_0^2q^{1/4}}
\notag
\end{equation}
and consequently
\begin{align}
2 C_0^2q^{1/4}G(1;q)&=\frac 1{\pi}\theta_4(0;\omega)\theta_1'(0;\omega)
\\
&=-\frac 1{\omega^2}\theta_4(0;-\frac 1{\omega})
\theta_3(0;-\frac 1{\omega})\theta_2(0;-\frac 1{\omega})^2
\notag
\end{align}
by Jacobi's transformation and the formula
\begin{equation}
\theta_1'(0;\omega)=\pi\theta_2(0;\omega)\theta_3(0;\omega)\theta_4(0;\omega).
\notag
\end{equation}
By Landen's transformation
\begin{equation}
\theta_4(0;-2/\omega)=\sqrt{\theta_3(0;-1/\omega)\theta_4(0;-1/\omega)}
\notag
\end{equation}
and hence
\begin{equation}
\frac{i\theta_4(0;-2/\omega)}{2\omega C_0^2q^{1/4}G(1;q)}=
-\frac{i\omega}{\theta_2(0;-1/\omega)^2
\sqrt{\theta_3(0;-1/\omega)\theta_4(0;-1/\omega)}}.
\notag
\end{equation}
If we insert this into (\ref{3.24}) we obtain (\ref{2.28}) and the theorem
is proved. That the leading behaviour of the kernel is given by (\ref{2.29})
follows from the series expansions of the $\theta$-functions:
\begin{align}\label{3.35}
\theta_1(x;\tau)&=i\sum_{n\in\mathbb{Z}}(-1)^ne^{\pi i\tau(n-1/2)^2+\pi i
(2n-1)x},
\\
\theta_2(x;\tau)&=\sum_{n\in\mathbb{Z}}e^{\pi i\tau(n-1/2)^2+\pi i
(2n-1)x},
\notag\\
\theta_3(x;\tau)&=\sum_{n\in\mathbb{Z}}e^{\pi i\tau n^2+2\pi inx},
\notag\\
\theta_4(x;\tau)&=\sum_{n\in\mathbb{Z}}(-1)^ne^{\pi i\tau n^2+2\pi inx}.
\notag\end{align}
\end{proof}

We should also prove theorem \ref{thm2.0}.

\begin{proof}
We will use the following criterion of Lenard, \cite{Len}, \cite{So}:

A family of localy integrable functions $\rho_k:I^k\to\mathbb{R}$,
$k=1,2,\dots$ are the correlation functions of some point process if and
only if the following two conditions are satisfied

\noindent
a) (Symmetry) For any $\sigma\in S_k$,
\begin{equation}
\rho_k(x_{\sigma(1)},\dots,x_{\sigma(k)})=\rho_k(x_1,\dots,x_k).
\notag
\end{equation}

\noindent
b) (Positivity) For any finite set of measurable bounded functions $\phi_k
:I^k\to\mathbb{R}$, $k=0,\dots, N$, with compact support, such that
\begin{equation}\label{3.36}
\phi_0+\sum_{k=1}^N\sum_{i_1\neq\dots\neq i_k} \phi_k(x_{i_1},\dots,x_{i_k})
\ge 0
\end{equation}
for all $(x_1,\dots,x_N)\in I^N$ it holds that
\begin{equation}\label{3.37}
\phi_0+\sum_{k=1}^N\int_{I^k}\phi_k(x_1,\dots,x_k)
\rho_k(x_1,\dots,x_k)dx_1\dots dx_k \ge 0.
\end{equation}
The uniform convergence of $K_N$ to $K$ on compact sets implies that $K$ is 
continuous and hence $\rho_k$ is locally integrable. It is also symmetric.
We know that (\ref{3.37}) holds with $\rho_{k,N}$ instead of $\rho_k$
since $\rho_{k,N}$ are the correlation functions of a point process.
Since all $\phi_k$ have compact support and are bounded we can take use the
uniform convergence of $K_N$ to $K$ and take $N\to\infty$ to get (\ref{3.37}).
This completes the proof.
\end{proof}

\section{Computation of the number variance}

We will first show how (\ref{2.38}) can be obtained from (\ref{2.31}).
By (\ref{2.35}) we want to compute
\begin{equation}\label{4.1}
\int_0^L dx\int_L^\infty dy f(x-y)+\int_0^L dx\int_{-\infty}^0 dy f(x-y),
\end{equation}
where
\begin{equation}\label{4.2}
f(x)=\frac{\sin^2\frac{\pi x}a}{\pi^2x^2}+\frac{d^2-(x/a)^2}{2\pi^2a^2
(d^2+(x/a)^2)^2}.
\end{equation}
Since $f(x)$ is even we see that (\ref{4.1}) equals
\begin{align}\label{4.3}
&2\int_{-L}^0dx\int_0^\infty dy f(x-y)=2\int_0^L\int_0^\infty dy f(x+y)
\\
&=2\int_0^L\left(\int_x^\infty f(y)dy\right)dx=
2\int_0^L xf(x)dx+2L\int_L^\infty f(x)dx.
\notag\end{align}
If we take $f(x)=(\sin^2\pi x)/\pi^2x^2$ in (\ref{4.3}) we get (\ref{2.34'})
using (\ref{4.14}) and (\ref{4.15}) below. Inserting (\ref{4.2}) 
into (\ref{4.3}) and computing the integrals we obtain (\ref{2.38}). A similar
computation using (\ref{2.31}) leads to (\ref{2.40}). To prove (\ref{2.41})
we have to show that
\begin{equation}\label{4.4}
\lim_{d\to\infty}\frac 2{\log d}\int_0^\infty \frac {u^2}{\sinh^2u}
\frac{\sin^22\pi ud}{u}du=1.
\end{equation}
Using $\sinh^2 u\ge u^2$ we get
\begin{align}
&\int_0^\infty\frac{2u^2}{\sinh^2u}\frac{\sin^22\pi ud}{u}du\le
\int_0^1\frac{1-\cos2\pi ud}u du+\int_1^\infty\frac{du}{\sinh^2u}
\\
&=\int_0^{2\pi d}\frac{1-\cos x}x dx+\int_1^\infty\frac{du}{\sinh^2u}
\notag\\
&=\log(2\pi d)+\gamma-\Ci(2\pi d)+
\int_1^\infty\frac{du}{\sinh^2u},
\notag
\end{align}
and hence
\begin{equation}
\limsup_{d\to\infty}\frac 2{\log d}\int_0^\infty \frac {u^2}{\sinh^2u}
\frac{\sin^22\pi ud}{u}du\le 1.
\notag
\end{equation}
Given $\epsilon>0$ we can choose $\delta>0$ so that $|x^{-1}\sinh x-1|\le
\epsilon$ if $|x|\le\delta$. Thus
\begin{align}
&\int_0^\infty \frac {2u^2}{\sinh^2u}\frac{\sin^22\pi ud}{u}du\ge
\frac 2{(1+\epsilon)^2}\int_0^\delta\frac{\sin^2 2\pi ud}u du
\\
&=\frac 1{(1+\epsilon)^2}(\log(2\pi d\delta)+\gamma-\Ci (2\pi d\delta)),
\notag
\end{align}
which gives
\begin{equation}
\liminf_{d\to\infty}\frac 2{\log d}\int_0^\infty \frac {u^2}{\sinh^2u}
\frac{\sin^22\pi ud}{u}du\ge 1,
\notag
\end{equation}
and we have proved (\ref{4.4}). This completes the proof of theorem
\ref{thm2.7}.

\begin{proof} {\it of theorem \ref{thm2.6}}.  
We have
\begin{align}
\text{Var}_{K_S}(\#[R,R+L])&=
\int_0^{R+L}dx\int_{R+L}^\infty dy K(x,y)K(y,x)
\notag\\
&+\int_R^{R+L} dx\int_{-\infty}^R dyK(x,y)K(y,x).
\notag
\end{align}
Set $L^\ast(x,y)=\frac 1{\pi}L(\frac x{\pi},\frac y{\pi})$, and
\begin{equation}
v(A,\theta)=\int_0^A dx\int_0^\infty dy L^\ast (-x+\theta,y+\theta)L^\ast
(y+\theta,-x+\theta).
\notag
\end{equation}
A computation shows that
\begin{equation}\label{4.5}
\text{Var}_{K_S}(\#[R,R+L])=
v(\frac{\pi L}a,\theta+\phi)+v(\frac{\pi L}a,-\theta).
\end{equation}
Hence we have to compute $v(A,\theta)$. If we set
\begin{equation}\label{4.6}
G(n,m,k)=\int_0^A dx\int_0^\infty dy\frac{e^{-2i((n+k)x+(m+k)y)}}
{(\pi dn+i(y+x))(\pi dm+i(y+x))},
\end{equation}
then
\begin{equation}\label{4.7}
v(A,\theta)=\frac 1{4\pi^2}\sum_{m,n\in\mathbb{Z}} e^{-dn(n-1)}
e^{2i(n-m)\theta}\left[e^{-dm(m+1)}G(n,m,0)+e^{-dm(m-1)}G(n,m,-1)\right]
+\text{c.c.}
\end{equation}
If we neglect the terms that are exponentially small in $d$ we get
\begin{align}\label{4.8}
4\pi ^2v(A,\theta)&=-G(0,0,0)+G(0,0,-1)+G(1,1,-1)
+e^{2i\theta}[G(1,0,-1)-G(0,-1,0)
\\
&-G(1,0,0)]
+e^{-2i\theta}G(0,1,-1)-e^{4i\theta}G(1,-1,0)+c.c.
\notag
\end{align}
(Note that $G(0,0,0)$ and $G(0,0,-1)$ are not individually convergent
but have to be considered together.)

We will now outline how (\ref{4.8}) can be computed without giving all
details. Set
\begin{equation}\label{4.9}
H_1(A;n,m)=\int_0^A\frac{1-e^{-2inx}}{\pi dm+ix}dx,
\end{equation}
$A>0$, $n,m\in\mathbb{Z}$, $d>0$ and
\begin{equation}\label{4.10}
H_2(A;n,m)=\int_A^\infty\frac{e^{-2inx}}{\pi dm+ix}dx,
\end{equation}
$A\ge 0$, $d>0$, $n\neq 0$, $n,m\in\mathbb{Z}$. Some computation now gives
\begin{align}\label{4.11}
&G(n,m,k)=\frac{1-e^{-2i(n-m)A}}{2\pi id(n-m)^2}[H_2(A;m+k,m)-H_2(A;m+k,n)]
\\
&+\frac 1{2\pi id(n-m)^2}[H_1(A;m+k,n)-H_1(A;m+k,m)+H_1(A;n+k,m)
-H_1(A;n+k,n)]
\notag
\end{align}
if $n\neq m$,
\begin{align}\label{4.12}
&G(n,n,k)=1-e^{-2i(n+k)A}-\frac 12\log(A^2+\pi^2d^2)+i\arctan\frac{\pi dn}A
+\log(\pi d)-i\frac{\pi}2\sgn(n)
\\
&-2(A-\pi idn)(n+k)H_2(A;n+k,n)-2\pi idn(n+k)H_2(0;n+k,n)+iH_1(A;n+k,n)
\notag
\end{align}
if $n\neq 0$ and
\begin{equation}\label{4.13}
G(0,0,k)-G(0,0,0)=1-e^{-2ikA}-2kAH_2(A;k,0)+H_1(A;k,0).
\end{equation}
The next step is to express the functions $H_1$ and $H_2$ in terms of
$f(z)$ and $g(z)$ defined by (\ref{2.37}), and in terms of the sine and
cosine integrals:
\begin{equation}\label{4.14}
\Si(A)=\int_0^A\frac{\sin y}{y}dy\quad;\quad \frac{\pi}2-\Si(A)=\int_A^\infty
\frac{\sin y}{y}dy,
\end{equation}
\begin{equation}\label{4.15}
\int_0^A\frac{1-\cos y}{y}dy=\gamma+\log A-\Ci(A)\quad;\quad 
\Ci(A)=-\int_A^\infty\frac{\cos y}{y}dy.
\end{equation}
After some computation we obtain
\begin{align}\label{4.16}
&H_1(A;n,m)=-i\log\sqrt{A^2+\pi^2d^2m^2}+i\log|\pi dm|
+\arg(A-\pi idm)
\\
&+\arg(\pi idm)+iG(-2\pi idm|n|)
-i(\cos 2nA)g(2A|n|-2\pi idm|n|)
\notag\\
&+i(\sin 2|n|A)f(2A|n|-2\pi idm|n|)
+f(-2\pi idm|n|)\sgn(n)
\notag\\
&-(\cos 2nA)f(2A|n|-2\pi idm|n|)\sgn(n)
-(\sin 2nA)g(2A|n|-2\pi idm|n|),
\notag\end{align}
if $n\neq 0$, $m\neq 0$,
\begin{equation}\label{4.17}
H_1(A;n,0)=-i(\gamma+\log(2A|n|)-\Ci(2|n|A)-\sgn(n)\Si(2A|n|),
\end{equation}
\begin{align}\label{4.18}
&H_2(A;n,m)=-i(\cos 2nA)g(2A|n|-2\pi idm|n|)+
i(\sin 2|n|A)f(2A|n|-2\pi idm|n|)
\\
&-(\cos 2nA)gf(2A|n|-2\pi idm|n|)\sgn(n)
-(\sin 2nA)g(2A|n|-2\pi idm|n|),
\notag\end{align}
if $m\neq 0$, $n\neq 0$,
\begin{equation}\label{4.19}
H_2(A;n,0)=i\Ci(2A|n|)-\sgn(n)(\frac{\pi}2-\Si(2A|n|))
\end{equation}
if $n\neq 0$, and finally
\begin{equation}\label{4.20}
H_2(A;0,0)-H_2(A;0,1)=i\log A-i\log\sqrt{A^2+\pi^2 d^2}+\arg(A-id).
\end{equation}
If we use these formulas in (\ref{4.8}), (\ref{4.11}), (\ref{4.12})
and (\ref{4.13}) we end up with (\ref{2.36}). The asymptotics for $f(z)$ 
and $g(z)$ are easy to obtain using integration by parts.
\end{proof}

\section{Proof of the approximation theorem}

\subsection{Main part of proof}

We will use the formulas (\ref{2.19})-(\ref{2.21}) for $K_S^{\text{ab}}(u,v)$.
Given a sequence $\{c_j\}_{j\ge 1}$ of complex numbers $\neq 0$, we
define the counting function,
\begin{equation}
n_c(t)=\#\{j\ge 1\,;\,|c_j|\le t\}.
\notag
\end{equation}
If $n_c(t)\le Ct^{1+\delta}$ for some $\delta$, $0\le\delta <1$, we can 
define the convergent canonical product
\begin{equation}\label{5.1}
P_c(z)=\prod_{j=1}^\infty \left(1-\frac z{c_j}\right) e^{z/c_j}.
\end{equation}
It follows from lemma \ref{lem5.1} below that $\xi_m$, as defined by 
(\ref{2.43}) is finite. Define
\begin{equation}\label{5.2}
c_j=c_{j,m}=y_{j+m}-y_m\quad;\quad b_j=b_{j,m}=y_m-y_{m-j}.
\end{equation}
We will first show that
\begin{equation}\label{5.3}
\frac{(w+y_m)F(w+y_m)}{(z+y_m)F(z+y_m)}=e^{\xi_m(w-z)}\frac {wP_c(w)P_b(-w)}
{zP_c(z)P_b(-z)},
\end{equation}
provided $z$ is not a zero of $P_c$ or $-z$ a zero of $P_b$. Note that
$P_b$ and $P_c$ are well defined by assumption (i) in the theorem. The
left hand side of (\ref{5.3}) is
\begin{align}
&\lim_{N\to\infty}\frac{w+y_m}{z+y_m}\prod_{j=1}^N\frac{y_j^2-(w+y_m)^2}
{y_j^2-(z+y_m)^2}
\\
&=\lim_{N\to\infty}\frac{w+y_m}{z+y_m}\prod_{j=1}^N
\frac{(y_j-(y_m+w))(y_j+y_m+w)}{(y_j-(y_m+z))(y_j+y_m+z)}
\notag\\
&=\lim_{N\to\infty}\frac{w+y_m}{z+y_m}\prod_{j=-N}^N
\frac{y_j-y_m-w}{y_j-y_m-z}=
\lim_{N\to\infty}\frac wz\prod_{j=-N}^{m-1}\frac{y_j-y_m-w}{y_j-y_m-z}
\prod_{j=m+1}^{N}\frac{y_j-y_m-w}{y_j-y_m-z}
\notag\\
&=\lim_{N\to\infty}\frac wz\prod_{j=1}^{N+m}\frac{y_m-y_{m-j}+w}
{y_m-y_{m-j}+z}\prod_{j=1}^{N-m}\frac{y_{j+m}-y_m-w}{y_{j+m}-y_m-z}
\notag\\
&=\lim_{N\to\infty}\frac wz\prod_{j=1}^{N+m}\frac{1+w/b_j}{1+z/b_j}
\prod_{j=1}^{N-m}\frac{1-w/c_j}{1-z/c_j}
\notag\\
&=\lim_{N\to\infty}\frac wz\prod_{j=1}^{N+m}
\frac{e^{-w(\frac 1{c_j}-\frac 1{b_j})}}{e^{-z(\frac 1{c_j}-\frac 1{b_j})}}
\prod_{j=N-m+1}^{N+m} e^{\frac 1{c_j}(z-w)}\frac{P_b(-w)P_c(w)}
{P_b(-z)P_c(z)},
\notag\end{align}
which gives the right hand side of (\ref{5.3}) since, $\sum_{j=N-m+1}^{N+m}
\frac 1{c_j}\to 0$ as $N\to\infty$ for a fixed $m$. 

Write $u_m=u-\zeta_m$, $v_m=v-\zeta_m$ (we write just $m$ instead of
$m(\alpha)$). Then
\begin{equation}\label{5.4}
|u_m|\le T+\lambda_m\quad,\quad |v_m|\le T+\lambda_m.
\end{equation}
Note that
\begin{equation}\label{5.5}
\frac{(w+y_m-v)^2}{2S}-\frac{(z+y_m-u)^2}{2S}=
\frac{(w-v_m)^2}{2S}-\frac{(z-u_m)^2}{2S}+\xi_m(w-z+u-v).
\end{equation}
If we make the change of variables $z\to z+y_m+u_m$, $w\to w+y_m+v_m$ in 
(\ref{2.19}) we obtain, using (\ref{5.3}),
\begin{align}\label{5.6}
K^\ast_{S,1}(u,v)&=\frac{e^{\xi_m(u-v)}}{(2\pi i)^2S}\int_{\Gamma_\sigma}
dw\int_{\gamma_M} e^{w^2/2S-z^2/2S}
\\
&\times\frac 1{z-w+u-v}
\frac{(w+v_m)P_c(w+v_m)P_b(-w-v_m)}{(z+u_m)P_c(z+u_m)P_b(-z-u_m)},
\notag\end{align}
where we have taken $L=\sigma+y_m+v_m$. The number 
$\sigma$ will be specified later
and satisfies $|\sigma|\le\lambda_m$. We will choose $M=\pi S/\lambda_m$ for
reasons that will be clear below. Note that $M\ge 1$ if $\alpha$ is large 
enough by our assumption on the allowed values of $S$. 
Below we will need the following estimates of the canonical products. 
Fix an $\epsilon>0$. There are constants $c_1,c_2>0$ such that
\begin{equation}\label{5.7}
|P_b(w)|,|P_c(w)|\le c_1e^{c_2(\lambda_m^{-1}|w|+|w|^{1+\delta})},
\end{equation}
for all $w\in\mathbb{C}$ and
\begin{equation}\label{5.8}
|P_b(z)|,|P_c(z)|\ge c_1^{-1}e^{-c_2(\lambda_m^{-1}|z|^{1+\epsilon}
+|z|^{1+\delta+\epsilon})},
\end{equation}
if $|\im z|\ge 1$.
These estimates are proved using the estimate
\begin{equation}
n_b(t)\le n_c(t)\le \lambda^{-1}t+Ct^{1+\delta},
\notag
\end{equation}
see (\ref{6.2}) below, and the following inequalities in \cite{Boas},
p. 19-22. If $x=\{x_k\}_{k=1}^\infty$ satisfies $n_x(t)\le Ct^{1+\delta}$,
$0\le\delta<1$, $t\ge 0$ and $r=|z|$, then
\begin{equation}\label{5.9}
\log|P_x(z)|\le 8\left\{r\int_0^r\frac{n_x(t)}{t^2}dt+r^2\int_r^\infty
\frac{n_x(t)}{t^3}dt\right\}
\end{equation}
for all $z\in\mathbb{C}$, and
\begin{equation}\label{5.10}
\log|P_x(z)|\ge -n_x(2r)\log(2r)+
\int_0^{2r}\frac{n_x(t)}{t}dt-8r^2\int_{2r}^\infty
\frac{n_x(t)}{t^3}dt-2r\int_0^{2r}\frac{n_x(t)}{t^2}dt,
\end{equation}
provided $|z-x_j|\ge 1$ for all $j\ge 1$.

Next we will prove an estimate which allows us to restrict the domain of
integration in (\ref{5.6}). Fix $R>0$. Introduce the following contours:
\begin{align}
&\Gamma_{\sigma,R}\,:\,[-R,R]\ni t\to\sigma+it,
\\
&\Gamma_{\sigma,R}^c\,:\,\mathbb{R}\setminus [-R,R]\ni t\to\sigma+it,
\notag\\
&\gamma_{M,R}\,:\,[-R,R]\ni t\to\mp t\pm iM,
\notag\\
&\gamma_{M,R}^c\,:\,\mathbb{R}\setminus [-R,R]\ni t\to\mp t\pm iM.
\notag
\end{align}
Let $(\gamma,\gamma')$ denote either $(\Gamma_\sigma,\gamma_{M,R}^c)$ or
$(\Gamma_{\sigma,R}^c,\gamma_M)$. We will show that there is a constant
$c_0$ such that if
\begin{equation}\label{5.11}
R\ge c_0\max\left\{\left(\frac{S}{\lambda_m}\right)^{1/(1-\epsilon)},
(T^{1+\delta+\epsilon}S)^{1/2},\left(\frac{T^{1+\epsilon}S}{\lambda_m}\right)
^{1/2}\right\}
\end{equation}
then
\begin{align}\label{5.12}
&\frac 1{S}\int_{\gamma} |dw|\int_{\gamma'}|dz|\left|e^{(w^2-z^2)/2S}\right|
\frac 1{|z-w+u-v|}\frac{|w+v_m|}{|z+u_m|}
\\
&\times\left|\frac{P_c(w+v_m)P_b(-w-v_m)}
{P_c(z+u_m)P_b(-z-u_m)}\right|\le \frac{C\lambda_m}{S^{3/2}} e^{-R^2/8S}.
\notag\end{align}
Write $g(t)=\lambda_m^{-1}t+t^{1+\delta}$ and 
$h(t)=\lambda_m^{-1}t^{1+\epsilon}+t^{1+\delta+\epsilon}$, and $I_R=[-R,R]$.
It follows from (\ref{5.7}) and (\ref{5.8}) that the integral in (\ref{5.12})
is
\begin{align}\label{5.13}
\le&\frac{C}{MS}\int_{\mathbb{R},I_R} ds\int_{I_R^c,\mathbb{R}} dt
e^{(-s^2+\sigma^2-t^2+M^2)/2S}\frac 1{\sqrt{(\mp t-\sigma+u-v)^2+(\mp s+M)^2}}
\\
&\times e^{cg(\sqrt{(\sigma+v_m)^2+s^2})+ch(\sqrt{(\mp t+u_m)^2+M^2})},
\notag\end{align}
for some constant $c>0$. Here we have used $|z+u_m|\ge M$ and 
$t\le \exp(g(t))$. We can now use $g(t+s)\le 2(g(t)+g(s))$ and similarly
for $h$ to see that the integral in (\ref{5.13}) is
\begin{align}
&\le\frac{C}{MS}\exp\left(\frac{\sigma^2+M^2}{2S} +C(g(|\sigma|)+g(|v_m|)
+h(|u_m|)+h(M))\right)
\\
&\times\int_{\mathbb{R},I_R^c}ds\int_{I_R^c,\mathbb{R}}dt\frac
{\exp\left(-\frac{s^2+t^2}{2S}+C(g(|s|)+h(|t|))\right)}
{\sqrt{(t-\sigma+u-v)^2+(s-M)^2}}
\notag\\
&\le\frac C{MS}\exp\left(\frac{\sigma^2+M^2}{2S}-\frac{R^2}{4S}+
C(g(|\sigma|)+g(|v_m|)+h(|u_m|)+h(M))\right)
\notag\\
&\times\int_{\mathbb{R}}ds\int_{\mathbb{R}}dt\frac
{\exp\left(-\frac{s^2+t^2}{4S}+C(g(|s|)+h(|t|))\right)}
{\sqrt{(t-\sigma+u-v)^2+(s-M)^2}}
\notag\\
&\le\frac{C}{MS^{1/2}}\exp\left(\frac{\sigma^2+M^2}{2S}-\frac{R^2}{4S}+
C(g(|\sigma|)+g(|v_m|)+h(|u_m|)+h(M))+\Delta\right),
\notag\end{align}
where
\begin{equation}
\Delta=\max_{(s,t)\in\mathbb{R}^2}\left(-\frac{s^2+t^2}4
+C(g(|s|\sqrt{S})+h(|t|\sqrt{S}))\right)
\le CS^{\frac{1+\epsilon}{1-\epsilon}}\lambda_m^{-\frac 2{1-\epsilon}},
\notag
\end{equation}
since $1\ge S\ge\lambda_m$ (essentially). We see that we need
\begin{equation}
\frac {R^2}{8S}\ge\frac{\sigma^2+M^2}{2S}
+C(g(|\sigma|)+g(|v_m|)+h(|u_m|)+h(M))+
 CS^{\frac{1+\epsilon}{1-\epsilon}}\lambda_m^{-\frac 2{1-\epsilon}},
\notag
\end{equation}
which holds if $R$ satisfies (\ref{5.11}). Here we have used 
$|\sigma|\le\lambda_m$, (\ref{5.4}) and $M=\pi S/\lambda_m$. This proves
(\ref{5.12}).

We will also need the following estimate. There is a constant $C$ such that
\begin{equation}\label{5.15}
\frac 1S
\int_{\Gamma_\sigma}|dw|\int_{\gamma_M} |dz|\left|e^{\frac{w^2-z^2}{2S}}
\right|\frac 1{|z-w+u-v|}\left|\frac{\sin\frac{\pi (w+v_m)}{\lambda_m}}
{\sin\frac{\pi (z+u_m)}{\lambda_m}}\right|\le \frac CS.
\end{equation}
The left hand side of (\ref{5.15}) is
\begin{equation}\label{5.16}
\frac 1S \int_{\mathbb{R}}d\tau\int_{\mathbb{R}} dt
\frac{e^{(-\tau^2-t^2+\sigma^2+M^2)/2S}}{|\mp t+u-v+i(\sigma+\tau\mp M)|}
\left|\frac{\sin\frac{\pi}{\lambda_m}(v_m+\sigma+i\tau)}
{\sin\frac{\pi}{\lambda_m}(u_m\pm t\mp M)}\right|.
\end{equation}
Now,
\begin{equation}
\left|\sin\frac{\pi}{\lambda_m}(v_m+\sigma+i\tau)\right|^2\le 2\cosh
\frac{2\pi\tau}{\lambda_m}\le 4e^{\frac{2\pi\tau}{\lambda_m}}
\end{equation}
and hence
\begin{equation}
\left|\sin\frac{\pi}{\lambda_m}(v_m+\sigma+i\tau)\right|\le 2
e^{\frac{\pi\tau}{\lambda_m}}.
\notag
\end{equation}
Also, by our choice of $M$,
\begin{equation}\label{5.17}
\left|\sin\frac{\pi}{\lambda_m}(u_m\pm t\mp M)\right|\ge\sinh\frac{\pi M}
{\lambda_m}=\sinh\frac{\pi^2S}{\lambda_m^2}
\end{equation}
Hence the integral in (\ref{5.16}) is
\begin{equation}
\le \frac CS
\frac{e^{\frac{\pi^2S}{\lambda_m^2}}}{\sinh\frac{\pi^2S}{\lambda_m^2}}
\int_{\mathbb{R}^2}\frac{e^{-\frac 1{2S}(\tau-\frac{\pi S}{\lambda_m})^2-
\frac {t^2}{2S}}}{\sqrt{(t-(u-v))^2+(\tau\pm\sigma\pm\frac{\pi S}{\lambda_m}
)^2}} d\tau dt.
\end{equation}
The contribution to the integral from $(t-(u-v))^2+(\tau\pm\sigma\pm\pi S/
\lambda_m)^2\le 1$ is $\le C$ and from the the contribution from the 
complementary region is $\le CS$. This proves (\ref{5.15}).

It follows from (\ref{5.12}) that
\begin{align}\label{5.18}
&e^{-\xi_m(u-v)}K^\ast_{S,1}(\underline{y};u,v)=
\frac 1{(2\pi i)^2S}\int_{\Gamma_{\sigma,R}} dw
\int_{\gamma_{M,R}} dze^{\frac{w^2-z^2}{2S}}\frac 1{z-w+u-v}
\\
&\times\frac{(w+v_m)P_c(w+v_m)P_b(-w-v_m)}{(z+u_m)P_c(z+u_m)P_b(-z-u_m)}
+\mathcal{R}_1,
\notag\end{align}
where
\begin{equation}\label{5.19}
|\mathcal{R}_1|\le\frac{C\lambda_m}{S^{3/2}}e^{-R^2/8S}
\end{equation}
provided $R$ satisfies (\ref{5.11}). Let $a_k=\lambda_m k$, $k\ge 1$. Then
\begin{equation}\label{5.20}
zP_a(z)P_a(-z)=\sin\frac{\pi z}{\lambda_m},
\end{equation}
and we write
\begin{align}\label{5.21}
&\frac{(w+v_m)P_c(w+v_m)P_b(-w-v_m)}{(z+u_m)P_c(z+u_m)P_b(-z-u_m)}=
\frac{\sin\frac{\pi}{\lambda_m}(w+v_m)}{\sin\frac{\pi}{\lambda_m}(z+u_m)}
\\
&+\left(\frac{P_c(w+v_m)}{P_a(w+v_m)}\frac{P_b(-w-v_m)}{P_a(-w-v_m)}
\frac{P_a(z+u_m)}{P_c(z+u_m)}\frac{P_a(-z-u_m)}{P_b(-z-u_m)}-1\right)
\frac{\sin\frac{\pi}{\lambda_m}(w+v_m)}{\sin\frac{\pi}{\lambda_m}(z+u_m)}
\notag
\end{align}
The argument above with $F(t)=\lambda_m^{-1}t$ gives
\begin{equation}\label{5.22}
K^\ast_{S,1}(\underline{\tilde{y}};u,v)=\int_{\Gamma_{\sigma,R}} dw
\frac 1{(2\pi)^2S}\int_{\gamma_{M,R}} dze^{\frac{w^2-z^2}{2S}}\frac 1{z-w+u-v}
\frac{\sin\frac{\pi}{\lambda_m}(w+v_m)}
{\sin\frac{\pi}{\lambda_m}(z+u_m)}
+\mathcal{R}_2,
\end{equation}
where
\begin{equation}\label{5.23}
|\mathcal{R}_2|\le\frac{C\lambda_m}{S^{3/2}}e^{-R^2/8S}
\end{equation}
provided $R$ satisfies (\ref{5.11}) (with an appropriate constant $c_0$
that does not depend on $\lambda_m$). Here 
$K^\ast_{S,1}(\underline{\tilde{y}};u,v)$ is given by (\ref{2.19}) 
with $L=\sigma+v_m$.

Set
\begin{align}\label{5.24}
\mathcal{R}_3&=\frac 1{(2\pi i)^2S}\int_{\Gamma_{\sigma,R}} dw
\int_{\gamma_{M,R}} dze^{\frac{w^2-z^2}{2S}}\frac 1{z-w+u-v}
\frac{\sin\frac{\pi}{\lambda_m}(w+v_m)}
{\sin\frac{\pi}{\lambda_m}(z+u_m)}
\\
&\times
\left(\frac{P_c(w+v_m)}{P_a(w+v_m)}\frac{P_b(-w-v_m)}{P_a(-w-v_m)}
\frac{P_a(z+u_m)}{P_c(z+u_m)}\frac{P_a(-z-u_m)}{P_b(-z-u_m)}-1\right).
\notag
\end{align}
Then, by (\ref{5.18}), (\ref{5.21}) and (\ref{5.22})
\begin{equation}\label{5.25}
e^{-\xi_m(u-v)}K^\ast_{S,1}(\underline{y};u,v)=
K^\ast_{S,1}(\underline{\tilde{y}};u,v)+\mathcal{R}_1-\mathcal{R}_2
+\mathcal{R}_3.
\end{equation}
We need an estimate of $\mathcal{R}_3$. For this we need estimates of
\begin{equation}
\left|\frac{P_c(w+v_m)}{P_a(w+v_m)}-1\right|
\quad\text{and}\quad
\left|\frac{P_a(z+u_m)}{P_c(z+u_m)}-1\right|
\end{equation}
and the same expression with $b$ instead of $c$ and a change of sign.
We have the identity
\begin{equation}\label{5.26}
\frac{P_c(w+v_m)}{P_a(w+v_m)}=\exp\left[\int_0^\infty\frac{(w+v_m)^2}
{w+v_m-t}\frac{n_c(t)-n_a(t)}{t^2} dt\right].
\end{equation}
Hence
\begin{equation}\label{5.27}
\left|\frac{P_c(w+v_m)}{P_a(w+v_m)}-1\right|\le
\exp\left[\int_0^\infty\left|\frac{(w+v_m)^2}
{w+v_m-t}\right|\frac{|n_c(t)-n_a(t)|}{t^2} dt\right]-1.
\end{equation}
Here we can use lemma \ref{lem5.2} below with 
\begin{equation}
g(t)=\frac{(\sigma+v_m)^2+s^2}{\sqrt{(\sigma+v_m-t)^2+s^2}},
\end{equation}
where $w=\sigma+is$.
Let $k_m=[v_m/\lambda_m-1/2]$ and choose $\sigma=\sigma_m=\lambda_m
(k_m+1/2)-v_m$. Then $|\sigma_m|\le\lambda_m$ and we have, by (\ref{5.4}), 
$|\sigma_m+v_m|\le T+2\lambda_m$. By assumption $T\le T_0\le1/4\sqrt{\eta_m}$,
(\ref{2.44}), and it follows that $k_m\le K=[(2\lambda_m\sqrt{\eta_m})^{-1}]$
if $m$ (i.e. $\alpha$) is large enough. Using the notation of lemma 
\ref{lem5.2} we see that if $|t-\lambda_m k|\le\alpha_k$, then
\begin{equation}\label{5.28}
|\sigma_m+v_m-t|\ge\frac{\lambda_m}8 (|k-k_m|+1),
\end{equation}
$1\le k\le K$. Here we have used that $\alpha_k\le\lambda_m/4$ if 
$1\le k\le K$. If $t\ge K\lambda_m$, then
\begin{equation}\label{5.29}
|\sigma_m+v_m-t|\ge \frac t4
\end{equation}
if $m$ is large enough. It follows from (\ref{5.28}) and (\ref{5.29}) that
\begin{equation}\label{5.30}
g(t)\le 8\frac{(T+2\lambda_m)^2+s^2}{\lambda_m(|k-k_m|+1)}
\end{equation}
if $|t-\lambda_mk|<a_k$ for $1\le k\le K$, and
\begin{equation}\label{5.31}
g(t)\le 4\frac{(T+2\lambda_m)^2+s^2}{t}
\end{equation}
if $t\ge K\lambda_m$. We ca now use lemma \ref{lem5.2} to conclude that
\begin{align}\label{5.32}
\int_0^\infty g(t)\frac{|n_c(t)-n_a(t)|}{t^2}dt
&\le\left[(T+2\lambda_m)^2+R^2\right]\left(\eta_m\log (\frac 1{\lambda_m
\sqrt{\eta_m}})+(\lambda_m\eta_m)^{1-\delta}\right)
\\
&\le C\left[(T+2\lambda_m)^2+R^2\right]m^{-\frac{(1-\delta)^2}{1+\delta}},
\notag\end{align}
where we have used $|s|\le R$, $\lambda_m^{-1}\le Cm^\delta$, $\eta_m\le C 
m^{-\frac{1-\delta}{1+\delta}}$ and $\lambda_m\le C$. Hence, by
(\ref{5.27}),
\begin{equation}\label{5.33}
\left|\frac{P_c(w+v_m)}{P_a(w+v_m)}-1\right|\le
\exp\left(Cm^{-\frac{(1-\delta)^2}{1+\delta}}
\left[(T+2\lambda_m)^2+R^2\right]\right)-1.
\end{equation}
A very similar computation using lemma \ref{lem5.4} instead gives
\begin{equation}\label{5.34}
\left|\frac{P_b(-w-v_m)}{P_a(-w-v_m)}-1\right|\le
\exp\left(Cm^{-\frac{(1-\delta)^2}{1+\delta}}
\left[(T+2\lambda_m)^2+R^2\right]\right)-1.
\end{equation}
We also have the estimate
\begin{equation}\label{5.35}
\left|\frac{P_a(z+u_m)}{P_c(z+u_m)}-1\right|\le
\exp\left[\int_0^\infty g(t)
\frac{|n_c(t)-n_a(t)|}{t^2} dt\right]-1.
\end{equation}
where now
\begin{equation}\label{5.36}
g(t)=\left|\frac{(z+u_m)^2}{z+u_m-t}\right|\le
\frac{(-s+u_m)^2+M^2}{\sqrt{(-s+u_m-t)^2+M^2}},
\end{equation}
if $z$ belongs to the upper part of $\gamma_{M,R}$ (the other case is
completely analogous). We have $|s|\le R$ and $2(R+T+\lambda_m)<
1/2\sqrt{\eta_m}$ if $m$ is sufficiently large by our assumptions on 
$R$ and the fact that $T\le T_0$. If $t\ge 1/2\sqrt{\eta_m}$ we get
\begin{equation}\label{5.37}
g(t)\le C\frac{R^2+(T+\lambda_m)^2}t.
\end{equation}
If $0\le t\le 1/2\sqrt{\eta_m}$, it could happen that $-s+u_m$ is close
to $t$. Here we use a similar estimate as above,
\begin{equation}\label{5.38}
g(t)\le C\frac{R^2+(T+\lambda_m)^2}{\sqrt{\lambda_m^2|k-k_\ast|^2+M^2}}
\end{equation}
with an appropriate $k_\ast$ (depending on  $-s+u_m$).
Using the estimates (\ref{5.37}) and (\ref{5.38}) in (\ref{5.35}) and
lemma \ref{lem5.2} we again get
\begin{equation}\label{5.39}
\left|\frac{P_a(z+u_m)}{P_c(z+u_m)}-1\right|\le
\exp\left(Cm^{-\frac{(1-\delta)^2}{1+\delta}}
[(T+2\lambda_m)^2+R^2]\right)-1,
\end{equation}
and similarly, using lemma \ref{lem5.4} instead,
\begin{equation}\label{5.40}
\left|\frac{P_a(-z-u_m)}{P_b(-z-u_m)}-1\right|\le
\exp\left(Cm^{-\frac{(1-\delta)^2}{1+\delta}}
[(T+2\lambda_m)^2+R^2]\right)-1,
\end{equation}
By our assumptions on $R$ and $T$ we see that the expression in the exponent
in (\ref{5.40}) is bounded by a constant. Hence by (\ref{5.15}), (\ref{5.24}),
(\ref{5.33}), (\ref{5.34}), (\ref{5.39}) and (\ref{5.40}) we get 
\begin{equation}\label{5.41}
|\mathcal{R}_3|\le Cm^{-\frac{(1-\delta)^2}{1+\delta}}(T^2+R^2).
\end{equation}
From (\ref{5.25}), (\ref{5.19}), (\ref{5.23}) and (\ref{5.41}) it follows that
\begin{equation}\label{5.42}
\left|e^{-\xi_m(u-v)}K^\ast_{S,1}(\underline{y};u,v)-
K^\ast_{S,1}(\tilde{\underline{y}};u_m,v_m)\right|
\le\frac CS m^{-\frac{(1-\delta)^2}{1+\delta}}(T^2+R^2)
+\frac {C\lambda_m}{S^{3/2}}e^{-R^2/8S}.
\end{equation}
Here $K^\ast_{S,1}(\underline{y};u,v)$ is given by (\ref{2.19}) with
$L=\sigma+y_m+v_m$ and $K^\ast_{S,1}(\tilde{\underline{y}};u_m,v_m)$ is given 
by (\ref{2.19}) with $L=\sigma+v_m$. Thus,
\begin{align}
&e^{-\xi_m(u-v)}K^\ast_{S,2}(\underline{y};u,v)
=\frac 1{2\pi i}e^{-\xi_m(u-v)}\int_{\sigma+y_m+v_m-Mi}^
{\sigma+y_m+v_m+Mi}e^{\frac 1{2S}((w-v)^2-(w-u)^2)} dw
\\
&=\frac 1{2\pi i}\int_{\sigma+v_m-Mi}^
{\sigma+v_m+Mi}e^{\frac 1{2S}((w-v)^2-(w-u)^2)} dw
=K^\ast_{S,2}(\tilde{\underline{y}};u_m,v_m).
\notag
\end{align}
Hence (\ref{5.42}) also gives
\begin{equation}\label{5.43}
\left|e^{-\xi_m(u-v)}K^\ast_{S}(\underline{y};u,v)-
K^\ast_{S}(\tilde{\underline{y}};u_m,v_m)\right|
\le\frac CS m^{-\frac{(1-\delta)^2}{1+\delta}}(T^2+R^2)
+\frac {C\lambda_m}{S^{3/2}}e^{-R^2/8S}.
\end{equation}

Note that
\begin{align}
K^\ast_{S,1}(\underline{y};-u,v)&=\frac 1{(2\pi i)^2S}\int_{\Gamma_L} dw
\int_{\gamma_M} dze^{\frac 1{2S}((w-v)^2-(z+u)^2)}\frac 1{z-w}
\frac{wF(z)}{zF(z)}
\notag\\
&=\frac 1{(2\pi i)^2S}\int_{\Gamma_L} dw
\int_{\gamma_M} dze^{\frac 1{2S}((w-v)^2-(z-u)^2)}\frac 1{z+w}
\frac{wF(z)}{zF(z)}
\notag
\end{align}
since $F(-z)=F(z)$. We can carry out the same type of computation as above to
see that (\ref{5.43}) also holds for $K^{\text{ab}}_S(\underline{y},u,v)$.
Now $K^{\text{ab}}_S(\tilde{\underline{y}},u_m,v_m)$ is approximated by 
$K_S(u_m,v_M)$, given by $K_S(u,v)=a^{-1}L_S(a^{-1}u,a^{-1}v)$ and $L_S$ as
in (\ref{2.23}), with error $\le \frac C{S\alpha}$, which is smaller than
the error term we have in the theorem. This completes the proof of the
approximation theorem.

\subsection{Some lemmas}

In the proof above we need some facts about certain numbers defined
in the theorem. 

\begin{lemma}\label{lem5.1}
The number $\xi_m$ defined by (\ref{2.43}) is finite. Also if we set
$\zeta_m=y_m-\xi_m S$, there is, for sufficiently large $\alpha$, an
$m=m(\alpha)$ such that
\begin{equation}\label{6.1}
|\zeta_{m(\alpha)}-\alpha|\le\lambda_{m(\alpha)}.
\end{equation}
\end{lemma}

\begin{proof}
In the proof we will need some rather simple facts which we will prove
later. They are immediate consequences of our assumptions on $F$.

\noindent
(a) $F'(t+s)\le F'(t)+F'(s)$, for all $t,s\ge 0$.

\noindent
(b) $F^{-1}(t+s)\le F^{-1}(t)+F^{-1}(s)$, for all $t,s\ge 0$.

\noindent
(c) $tF'(t)\le 4F(t)$, for all $t\ge 0$.

\noindent
(d) $tF''(t)\le F'(t)$, for all $t\ge 0$.

\noindent
(e) $|F(F^{-1}(m)+t)-m-\lambda_m^{-1}|\le\eta_m t^2$ for all $m\ge 1$, 
$t\ge 0$.

\noindent
(e) $|F(F^{-1}(m)+t)-m|\le \lambda_m^{-1} t+Ct^{1+\delta}$ for all $m\ge 1$, 
$t\ge 0$.

\noindent
(f) $\lambda_m-\lambda_{m+j}\le\eta_m\lambda_M^3 j$, $j\ge 1$

\noindent
(g) $\lambda_{m+1}\le y_{m+1}-y_m\le\lambda_m$.

Let $c_{j,m}$ and $b_{j,m}$ be defined by (\ref{5.2}). Then
$n_c(t)=[F(F^{-1}(m)+t)-m]$, where $[\cdot]$ denotes the
integer part, and hence by (f),
\begin{equation}\label{6.2}
n_c(t)\le\lambda^{-1}t+Ct^{1+\delta}
\end{equation}
for $t\ge 0$. We have
\begin{equation}\label{6.3}
\xi_m=\sum_{j=1}^\infty\left(\frac 1{c_{j,m}}-\frac 1{b_{j,m}}\right)
\end{equation}
and since $F$ is convex, $b_{j,m}\ge c_{j,m}$. Thus $0\le b_{j,m}-c_{j,m}
=2y_m+y_{j,m}-y_{j+m}\le 2y_m$ and we see from (\ref{6.3}) that $\xi_m\ge 0$
and
\begin{equation}\label{6.4}
\xi_m=\sum_{j=1}^\infty\frac{2y_m}{(y_{j+m}-y_m)(y_m+y_{j-m})}.
\end{equation}
Since $F(x)\le Cx^{1+\delta}$ we have 
\begin{equation}\label{6.4'}
y_j=F^{-1}(j)\ge Cj^{1/(1+\delta)}
\end{equation}
and it follows that the series in (\ref{6.4}) is convergent.

To prove the other statement in the lemma, (\ref{6.1}), we want to estimate 
$|\xi_m-\xi_{m+1}|$. From (\ref{6.3}) we have
\begin{equation}\label{6.5}
\xi_m-\xi_{m+1}=\sum_{j=1}^\infty
\left[\frac{\Delta y_{m+j}-\Delta y_m}{c_{j,m}c_{j,m+1}}+
\frac{\Delta y_{m-j}-\Delta y_m}{b_{j,m}b_{j,m+1}}\right],
\end{equation}
where we have used the notation $\Delta y_k=y_{k+1}-y_k$.
If we take $t=\Delta y_m$ in (e) and use (h) we get
\begin{equation}\label{6.6}
|\lambda_m-\Delta y_m|\le\lambda_m^3\eta_m
\end{equation}
and together with (g) this gives
\begin{equation}\label{6.7}
|\Delta y_{m+j}-\Delta y_m|\le\lambda_m^3\eta_m (j+2),
\end{equation}
for $j,m\ge 1$. Since $\lambda_m$ is decreasing in $m$, (h) gives
\begin{equation}\label{6.8}
c_{j,m}=y_{m+j}-y_m\ge j\lambda_{m+j}.
\end{equation}
Hence, if $1\le j\le m$,
\begin{equation}\label{6.9}
c_{j,m}c_{j,m+1}\ge j^2\lambda_{2m}^2\ge\frac 14 j^2\lambda_m^2.
\end{equation}
Here we have used
\begin{equation}\label{6.9'}
\frac{\lambda_m}{\lambda_{2m}}=\frac{F'(F^{-1}(2m))}{F'(F^{-1}(m))}
\le\frac{F'(F^{-1}(m)+F^{-1}(m))}{F'(F^{-1}(m))}\le 2
\end{equation}
by (a) and (b). Combining (\ref{6.7}) and (\ref{6.9}) we get
\begin{equation}\label{6.10}
\left|\sum_{j=1}^{m-1}\frac{\Delta y_{m+j}-\Delta y_m}{c_{j,m}c_{j,m+1}}
\right|\le 4\lambda_m\eta_m\sum_{j=1}^m\frac{j+2}{j^2}
\le 36\lambda_m\eta_m\log m
\end{equation}
if $m\ge 2$. From (c), (d) and (\ref{6.4'}) we get
\begin{equation}\label{6.11}
\eta_m=F''(y_m)\le\frac{4F(y_m)}{y_m^2}\le Cm^{-\frac{1-\delta}{1+\delta}},
\end{equation}
and thus (\ref{6.10}) gives
\begin{equation}\label{6.12}
\left|\sum_{j=1}^{m-1}\frac{\Delta y_{m+j}-\Delta y_m}{c_{j,m}c_{j,m+1}}
\right|\le C(\log m)m^{-\frac{1-\delta}{1+\delta}}\lambda_m.
\end{equation}
By (h) and the fact that $\lambda_m$ is dereasing we get
$|\Delta y_{m+j}-\Delta y_m|\le 2\lambda_m$. Hence, by (\ref{6.2}),
\begin{align}\label{6.13}
&\left|\sum_{j=m}^{\infty}\frac{\Delta y_{m+j}-\Delta y_m}{c_{j,m}c_{j,m+1}}
\right|\le 2\lambda_m\sum_{j=m}^{\infty}\frac 1{c_{j,m}^2}=
2\lambda_m\int_{c_{m,m}-}^\infty\frac{dn_c(t)}{t^2}
\\
&\le 4\lambda_m\int_{c_{m,m}}^\infty\frac{\lambda_m^{-1}+Ct^{1+\delta}}
{t^3} dt\le C\lambda_m\left(\frac 1{\lambda_mc_{m,m}}+\frac 1{c_{m,m}^
{1-\delta}}\right).
\notag\end{align}
By (h) and (\ref{6.9'}), $c_{m,m}=y_{2m}-y_m\ge\lambda_{2m}\ge\frac 12 
m\lambda_m$. It follows that the right hand side of (\ref{6.13}) is
$\le C\lambda_m((m\lambda_m^2)^{-1}+(m\lambda_m)^{\delta-1})$. Now, by (c),
\begin{equation}
\frac{F(t)}{F'(t)^2}=\frac 1{F(t)}\left(\frac{F(t)}{F'(t)}\right)^2\ge
\frac {t^2}{16F(t)}\ge Ct^{1-\delta}
\notag
\end{equation}
and consequently by (\ref{6.4'})
\begin{equation}\label{6.13'}
m\lambda_m^2=\frac{F(F^{-1}(m))}{F'(F^{-1}(m))^2}\ge Cy_m^{1-\delta}
\ge Cm^{\frac{1-\delta}{1+\delta}}.
\end{equation}
We find
\begin{equation}\label{6.14}
\left|\sum_{j=m}^{\infty}\frac{\Delta y_{m+j}-\Delta y_m}{c_{j,m}c_{j,m+1}}
\right|\le  Cm^{-\frac{1-\delta}{1+\delta}}\lambda_m.
\end{equation}
Combining (\ref{6.12}) and (\ref{6.14}) we find
\begin{equation}\label{6.15}
\left|\sum_{j=1}^{\infty}\frac{\Delta y_{m+j}-\Delta y_m}{c_{j,m}c_{j,m+1}}
\right|\le  Cm^{-\frac{1-\delta}{1+\delta}}\lambda_m.
\end{equation}

The same argument that led to (\ref{6.12}) gives
\begin{equation}\label{6.16}
\left|\sum_{j=1}^{m/2}\frac{\Delta y_{m-j}-\Delta y_m}{b_{j,m}b_{j,m+1}}
\right|\le  Cm^{-\frac{1-\delta}{1+\delta}}\lambda_m.
\end{equation}
Using $b_{j,m}b_{j,m+1}\ge c_{j,m}c_{j,m+1}$ and $\Delta y_{m-j}=
\Delta y_{j-m-1}$ for $j>m$, we see by similar arguments as above that
\begin{equation}\label{6.17}
\left|\sum_{j=3m/2}^{\infty}\frac{\Delta y_{m-j}-\Delta y_m}{b_{j,m}b_{j,m+1}}
\right|+\sum_{j=m/2}^{\infty}\frac{\Delta y_m}{b_{j,m}b_{j,m+1}}
\le  Cm^{-\frac{1-\delta}{1+\delta}}\lambda_m.
\end{equation}

It remains to consider
\begin{equation}\label{6.18}
\sum_{j=m/2}^{3m/2}\frac{|\Delta y_{m-j}|}{b_{j,m}b_{j,m+1}}
\le 2\sum_{j=1}^{m/2}\frac{\Delta y_j}{b_{m-j,m}b_{m-j,m+1}}.
\end{equation}
Now, by (h), $b_{m-j,m}=y_m-y_j\ge (m-j)\lambda_m$ and hence, by 
(\ref{6.13'}), (c)
and straightforward estimates
\begin{equation}
\sum_{j=1}^{m/2}\frac{\Delta y_j}{b_{m-j,m}b_{m-j,m+1}}
\le C\frac{y_m}{m\lambda_m}m^{-\frac{1-\delta}{1+\delta}}\lambda_m
\le Cm^{-\frac{1-\delta}{1+\delta}}\lambda_m.
\end{equation}
Combining this with (\ref{6.16}), (\ref{6.17}) and (\ref{6.18}) we obtain
\begin{equation}\label{6.19}
\left|\sum_{j=1}^{\infty}\frac{\Delta y_{m-j}-\Delta y_m}{b_{j,m}b_{j,m+1}}
\right| \le C(\log m)m^{-\frac{1-\delta}{1+\delta}}\lambda_m.
\end{equation}

From (\ref{6.5}), (\ref{6.15}) and (\ref{6.19}) we now get the disired 
estimate
\begin{equation}\label{6.20}
|\xi_m-\xi_{m+1}|\le 
C(\log m+1)m^{-\frac{1-\delta}{1+\delta}}\lambda_m
\end{equation}
for $m\ge 1$.
Using (\ref{6.20}) and (h) we get, for $1\le r\le m-1$,
\begin{align}
|\xi_m|&\le \xi_1+\sum_{j=1}^{m-1}|\xi_{k+1}-\xi_k|\le
\xi_1+C\sum_{j=1}^{m-1}(\log m+1)m^{-\frac{1-\delta}{1+\delta}}\Delta y_{k-1}
\\
&\le\xi_1+Cy_r+C(\log r+1)r^{-\frac{1-\delta}{1+\delta}}(y_{m-1}-y_{r-1}).
\notag\end{align}
Consequently $\xi_m/y_m\to 0$ as $m\to\infty$ and since $S\le 1$, we see that
$\zeta_m=y_m-S\xi_m\to\infty$ as $m\to\infty$. Since $\Delta y_m/\lambda_m
\to 1$ as $m\to\infty$ it follows from (\ref{6.20}) that
\begin{equation}
\left|\frac{\zeta_{m+1}-\zeta_m}{\lambda_m}\right|=
\left|\frac{\Delta y_m}{\lambda_m}-\frac{S(\xi_{m+1}-\xi_m)}{\lambda_m}\right|
\to 1
\notag
\end{equation}
as $m\to\infty$. Hence $|\zeta_{m+1}-\zeta_m|\le 3\lambda_m/2$ if $m$ is
sufficiently large. If we take $\alpha$ sufficiently large the closest $\xi_m$
is thus within distance $3\lambda_m/4$ or $3\lambda_{m+1}/4$, which is 
$\le\lambda_m$ and we take this $m$ as our
$m(\alpha)$.
\end{proof}

Our next lemma is

\begin{lemma}\label{lem5.2}
Let $a_j=\lambda_m j$, $j\ge 1$ and $c_j=c_{j,m}$, $j\ge 1$ as above.
If $g$ is a bounded measurable function on $[0,\infty)$ we have
\begin{align}\label{6.21}
&\left|\int_0^\infty g(t)\frac{|n_c(t)-n_a(t)|}{t^2} dt\right|
\le\sum_{k=1}^K\int_{k\lambda_m-\alpha_k}^{k\lambda_m+\alpha_k}
\frac{|g(t)|}{t^2}dt
\\
&+5\eta_m\int_{1/2\sqrt{\eta_m}}^{1/\lambda_m\eta_m} |g(t)|dt
+\frac 2{\lambda_m}\int_{1/\lambda_m\eta_m}^\infty\frac{|g(t)|}{t}dt
+C\int_{1/\lambda_m\eta_m}^\infty\frac{|g(t)|}{t^{1-\delta}} dt,
\notag\end{align}
where $\alpha_k=\lambda_m^3\eta_m(k+1)^2$, $K=[(2\lambda_m\sqrt{\eta_m}
)^{-1}]$.
\end{lemma}

\begin{proof} We have $n_c(t)=[F(t+F^{-1}(m))-m]$ and $n_a(t)=[\lambda_m
^{-1} t]$. The proof of (\ref{6.21}) is based on the following claim which 
we will prove below.

\begin{claim}\label{cl5.3}
Assume that $0\le t\le \lambda_mK$, $1\le k\le K$. Then $n_c(t)=n_a(t)$ if
$|t-\lambda_mk|\ge\alpha_k$ and $|n_c(t)-n_a(t)|\le 1$ if $|t-\lambda_mk|
\le \alpha_k$.
\end{claim}

Using the claim we have that the left hand side of (\ref{6.21}) is
\begin{equation}\label{6.22}
\le\sum_{k=1}^K\int_{k\lambda_m-\alpha_k}^{k\lambda_m+\alpha_k}
\frac{|g(t)|}{t^2}dt+
\int_{\lambda_mK}^\infty|g(t)|\frac{|n_c(t)-n_a(t)|}{t^2}dt.
\end{equation}
It follows from (e) that $|n_c(t)-n_a(t)|\le\eta_m t^2+1$ and from (f)
we get
\begin{equation}
|n_c(t)-n_a(t)|\le n_c(t)+n_a(t)\le \frac 2{\lambda_m} t+Ct^{1+\delta}.
\notag
\end{equation}
Using these estimates in the second integral in (\ref{6.22}) we obtain
(\ref{6.21}).
\end{proof}

\begin{proof} {\it of claim \ref{cl5.3}}.
If $0\le t\le\lambda_m K$, then $t\in [\lambda_mk,\lambda_m(k+1)]$ for some
$k$, $0\le k<K$. If furthermore $t\in[\lambda_m k+\alpha_k,
\lambda_m(k+1)-\alpha_{k+1}]$, then $[\lambda^{-1}t]=k$, since
$\lambda_m^{-1}\alpha_k=\lambda_m^2\eta_m(k+1)^2\le\lambda_m^2\eta_m
K^2\le 1/4$, by our choice of $K$. Also, $\lambda_m^{-1}t-[\lambda_m^{-1}t]
\ge\lambda_m^{-1}\alpha_k>\eta_m t^2$ since $t<\lambda_m(k+1)$, and
$\lambda_m^{-1}t-([\lambda_m^{-1}t]+1)\le-\lambda_m^{-1}\alpha_{k+1}<
-\eta_mt^2$ since $t<\lambda_m(k+1)$. Combined with (e) this gives
$0<F(t+F^{-1}(m))-m-[\lambda_m^{-1}t]<1$, i.e. $n_c(t)=n_a(t)$. On the other 
hand, if $0\le t\le\lambda_mK$ and $|t-\lambda_mk|<\alpha_k$, then $n_c(t)$
and $n_a(t)$ can differ by at most 1.
\end{proof}

We have a similar lemma for $n_b$ instead.

\begin{lemma}\label{lem5.4}
Using the same notation as in lemma \ref{lem5.2} we have
\begin{align}\label{6.23}
&\left|\int_0^\infty g(t)\frac{|n_b(t)-n_a(t)|}{t^2} dt\right|
\le\sum_{k=1}^K\int_{k\lambda_m-\alpha_k}^{k\lambda_m+\alpha_k}
\frac{|g(t)|}{t^2}dt
\\
&+\int_{1/4\sqrt{\eta_{m/2}}}^{m\lambda_m/2} |g(t)|(\frac 1{t^2}+2\eta_{m/2})
dt
+\frac 2{\lambda_m}\int_{m\lambda_m/2}^\infty\frac{|g(t)|}{t}dt
+C\int_{m\lambda_m/2}^\infty\frac{|g(t)|}{t^{1-\delta}} dt,
\notag\end{align}
where the second integral in the right hand side is present only if
$m\lambda_m/2>1/4\sqrt{\eta_{m/2}}$.
\end{lemma}

\begin{proof}
An integration by parts and the fact that $F''$ is decreasing gives
\begin{equation}
|m-F(F^{-1}(m)-t)-\lambda_m^{-1}t|\le 2\eta_{m/2}t^2
\notag
\end{equation}
if $0\le t\le y_m-y_{m/2}$. Note that by (h), $ y_m-y_{m/2}\ge m\lambda_m/2$.
Now, $n_b(t)$ equals $[m-F(F^{-1}(m)-t)]$ if $0\le t\le y_m$ and
$m+n_y(t-y_m)$ if $t>y_m$. It is clear that $n_b(t)\le\lambda_m^{-1} t$ for
$0\le t\le y_m$. If $t>y_m$, then using $n_y(t)\le F(t)$ we find
\begin{equation}
n_b(t)\le F(F^{-1}(m))+F(t-F^{-1}(m))\le F(t)\le Ct^{1+\delta}.
\notag
\end{equation}
The proof of (\ref{6.23}) now proceeds in the same way as the proof of
lemma \ref{lem5.2}.
\end{proof}

It remains to prove the statements in the begining of this subsection.
(a) and (b) are immediate consequences of our assumptions on $F$. To prove
(c) write
\begin{equation}
F(t)=\int_0^t F'(s)ds\ge \int_{t/2}^t F'(s)ds\ge \frac t2 F'(\frac t2),
\notag
\end{equation}
since $F'$ is increasing. By (a), $F'(t)\le 2F'(t/2)$, and (c) follows. The 
statement (d) follows from the fact that $F''$ is decreasing. To prove (e) 
write
\begin{align}
|F(F^{-1}(m)+t)-m-\lambda_m^{-1}t|&=\left|\int_{F^{-1}(m)}^{F^{-1}(m)+t}
(F^{-1}(m)+t-s)F''(s)ds\right|
\\
&\le F''(F^{-1}(m))t^2=\eta_mt^2,
\notag
\end{align}
since $F''$ is decreasing. For (f) write
\begin{align}
&F(F^{-1}(m)+t)-m=\int_{F^{-1}(m)}^{F^{-1}(m)+t}
F'(s)ds\le tF'(F^{-1}(m)+t)
\\
&\le tF'((F^{-1}(m))+tF'(t)\le\lambda_m^{-1} t+Ct^{1+\delta},
\notag
\end{align}
by (a), (c) and the fact that $F'$ is increasing. To prove (g) write
\begin{align}
&\lambda_m-\lambda_{m+j}=\frac{F'(F^{-1}(m+j))-F'(F^{-1}(m))}
{F'(F^{-1}(m))F'(F^{-1}(m+j))}
\le\frac 1{F'(y_m)^2}\int_{y_m}^{y_{m+j}}F''(s)ds
\\
&\le \eta_m\lambda_m^2(y_{m+j}-y_m)=
\eta_m\lambda_m^2\int_m^{m+j}\frac{dt}{F'(F^{-1}(t))}\le\eta_m\lambda_m^3 j.
\notag
\end{align}
Finally, to prove (h) we write
\begin{equation}
y_{m+1}-y_m=\int_m^{m+1}\frac{dt}{F'(F^{-1}(t))}.
\notag
\end{equation}
Since $F'\circ F^{-1}$ is increasing the right hand side is $\le\lambda_m$
and $\ge\lambda_{m+1}$.

\noindent
{\bf Acknowledgement}: 
This work was supported by the Swedish Science Research Council (VR)
and the G\"oran Gustafsson Foundation (KVA).

\end{document}